\documentclass[titlepage,11pt]{article}
\usepackage{latexsym}
\usepackage{amsfonts}
\usepackage{amsmath}
\usepackage{comment}
\usepackage{tikz}
\usepackage[T1]{fontenc}
\usepackage{lmodern}
\usetikzlibrary{decorations.pathreplacing,calligraphy}
\newcommand\blackslug{\hbox{\hskip 1pt \vrule width 4pt height 8pt depth 1.5pt
        \hskip 1pt}}
\newcommand\bbox{\hfill \quad \blackslug \bigbreak}

\renewcommand{\sc}{\operatorname{sc}}
\newcommand{\insc}{\operatorname{insc}}
\newcommand{\outsc}{\operatorname{outsc}}
\newcommand{\IP} {\operatorname{IP}}
\newcommand{\SP} {\operatorname{SP}}

\newcommand{\OP} {\operatorname{OP}}
\newcommand{\F} {\operatorname{FP}}

\newcommand{\mac}{\mathcal}

\def\LL{,\ldots,}
\def\cupcup{\cup\cdots\cup}

\title{Rough structure. I. Exclusion among series-parallel graphs}
\author{Maria Chudnovsky\thanks{Supported by NSF grants DMS-2348219 and   CCF-2505100,
AFOSR grant FA9550-25-1-0275, and a Guggenheim Fellowship.}\\
Princeton University, Princeton, NJ 08544
\\
\\
Julien Codsi\\
Princeton University, Princeton, NJ 08544
\\
\\
Alex Divoux\\
Princeton University, Princeton, NJ 08544
\\
\\
Paul Seymour\thanks{Supported by AFOSR grant
FA9550-22-1-0234, and NSF grant  DMS-2154169.}\\
Princeton University, Princeton, NJ 08544
\\
\\
Liana Yepremyan\thanks{Supported by NSF grant DMS-2247013.}\\
Emory University, Atlanta, GA 30322
}

\date{September 15, 2026; revised \today}

\newtheorem{thm}{}[section]

\newcommand{\Proof}{\noindent{\bf Proof.}\ \ }

\begin{document}
\maketitle
\begin{abstract}
This is the first in a series of papers about rough structure theorems: graph properties that are implied by the exclusion of some graph as a minor, and imply the exclusion of some larger graph of the same 
type. Here we look at the structure of the series-parallel graphs that do not contain a given series-parallel graph as a minor. We prove that for every series-parallel graph $H$, there exists $m$ such that 
every simple two-connected series-parallel graph that does not contain $H$ as a minor admits a tree-decomposition of height at most $m$, where each ``torso'' is either an outerplanar 
graph or a planar dual of one. We prove similar theorems about the series-parallel graphs that exclude an outerplanar graph, the series-parallel graphs that exclude a dual of an outerplanar graph, the graphs of tree-width three that exclude a series-parallel graph, and
the general graphs that exclude a fan.
\end{abstract}

\section{Introduction}
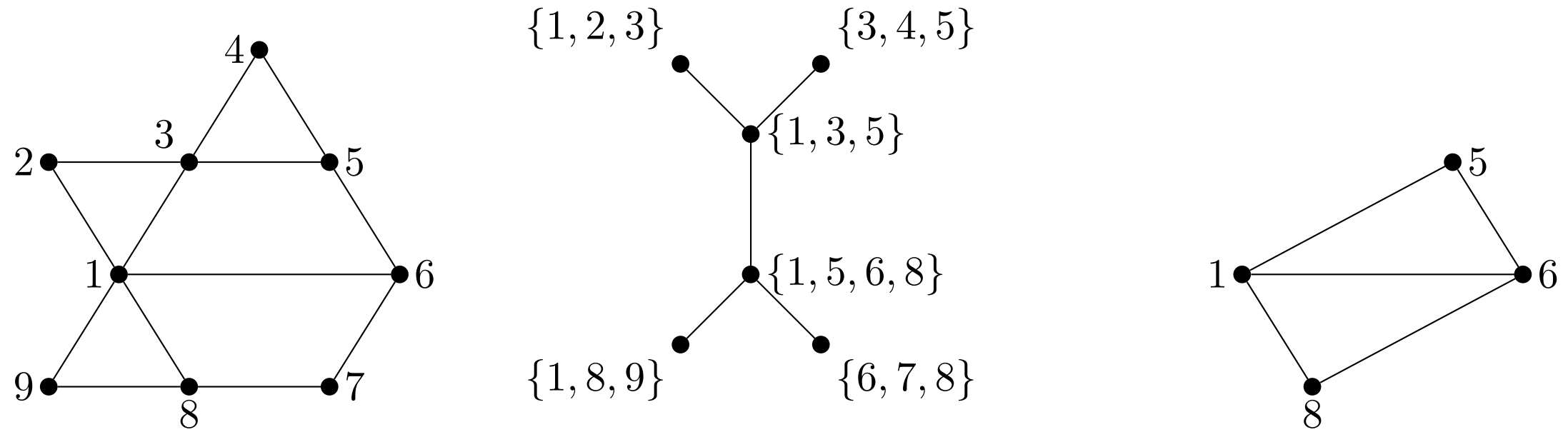
\begin{figure}[h]
\begin{center}
\begin{tikzpicture}[scale = 1.3]
\tikzstyle{every node}=[inner sep=1.5pt, fill=black,circle,draw]
\def\h{.8}
\node (1) at (-1,0) {};
\node (2) at (-1.5,\h) {};
\node (3) at (-.5,\h) {};
\node (4) at (0,2*\h) {};
\node (5) at (.5,\h) {};
\node (6) at (1,0) {};
\node (7) at (.5,-\h) {};
\node (8) at (-.5,-\h) {};
\node (9) at (-1.5,-\h) {};

\node (a) at (3,1.5) {};
\node (b) at (4,1.5) {};
\node (c) at (3.5,1) {};
\node (d) at (3.5,0) {};
\node (e) at (4,-.5) {};
\node (f) at (3,-.5) {};

\node (11) at (7,0) {};
\node (51) at (8.5,\h) {};
\node (61) at (9,0) {};
\node (81) at (7.5,-\h) {};

\foreach \from/\to in {1/2,1/3,2/3,3/4,3/5,4/5,5/6,1/6,6/7,7/8,1/8,1/9,8/9, a/c,b/c,c/d,d/e,d/f, 11/51,11/61,11/81,51/61,61/81}
\draw [-] (\from) -- (\to);

\tikzstyle{every node}=[]
\draw (1) node[left] {$1$};
\draw (2) node[left] {$2$};
\draw (3) node[above left] {$3$};
\draw (4) node[left] {$4$};
\draw (5) node[right] {$5$};
\draw (6) node[right] {$6$};
\draw (7) node[right] {$7$};
\draw (8) node[below] {$8$};
\draw (9) node[left] {$9$};

\draw (a) node[above left] {$\{1,2,3\}$};
\draw (b) node[above right] {$\{3,4,5\}$};
\draw (c) node[right] {$\{1,3,5\}$};
\draw (d) node[right] {$\{1,5,6,8\}$};
\draw (e) node[below right] {$\{6,7,8\}$};
\draw (f) node[below left] {$\{1,8,9\}$};

\draw (11) node[left] {$1$};
\draw (51) node[right] {$5$};
\draw (61) node[right] {$6$};
\draw (81) node[below] {$8$};

\end{tikzpicture}
\end{center}
\caption{A graph and a tree-decomposition of it. The right-hand graph is one of the torsos.} \label{fig:torso}
\end{figure}

All graphs in this paper are finite, and may have loops or parallel edges.
A {\em tree-decomposition} $\mac T$ of a graph $G$ is a pair $(T,(W_t:t\in V(T)))$, where $T$ is a tree, and $W_t\subseteq V(G)$ for each $t\in V(T)$, such that 
\begin{itemize}
\item $G=\bigcup_{t\in V(T)}G[W_t]$, and
\item if $t,t',t''\in V(T)$ and $t'$ lies on the path between $t, t''$, then $W_t\cap W_{t''}\subseteq W_{t'}$.
\end{itemize}
(See Figure \ref{fig:torso}.)

The {\em width} of $\mac T$ is the maximum of $|W_t|-1$ over all $t\in V(T)$, and the {\em tree-width} of $G$ is the minimum width of its 
tree-decompositions. The simple graphs of tree-width at most one are the forests, and the graphs of tree-width at most two are the series-parallel graphs.
There is a well-known theorem of Robertson and Seymour~\cite{GM5}:
\begin{thm}\label{GM5}
For every planar graph $H$,
there exists $\ell$ such that every graph not containing $H$ as a minor has tree-width at most $\ell$.
\end{thm}
(A graph is a {\em minor} of another if the first can be obtained from a subgraph of the second by edge-contraction.)
See~\cite{diestel} for a simpler proof, and~\cite{chuzhoy} (and, very recently,~\cite{haj}) for the best dependence between $H$ and $\ell$ known. 

This is best possible in a sense, because for every non-planar graph $H$ there is no such $\ell$ (to see this, observe that large planar grids have large tree-width and do not contain $H$
as a minor if $H$ is non-planar). And it has a converse in a sense; because for all $\ell$, there is a planar graph $H$ such that every graph with tree-width at most $\ell$ does not contain $H$ as a minor.
(To see this, choose a grid with tree-width more than $\ell$, and observe that every minor of a graph with tree-width at most $\ell$ also has tree-width at most $\ell$.)

There is a similar theorem about excluding a forest,
also due to Robertson and Seymour~\cite{GM1} (a graph has {\em path-width} at most $k$ if it admits a tree-decomposition $(T,(W_t:t\in V(T)))$ of width at most $k$ where $T$ is a path):
\begin{thm}\label{GM1}
For every forest $H$
there exists $\ell$ such that every graph not containing $H$ as a minor 
has path-width at most $\ell$.
\end{thm}
(See Diestel~\cite{diestel2} for a short proof.)
Again, this is ``best possible'', and has a ``converse'', in the same senses as before, since uniform binary trees of large depth have large path-width.
Such theorems are called ``rough structure theorems''.

Let us say an {\em ideal} of graphs is a class of graphs closed under isomorphism and under taking minors: thus, the class of all planar graphs is an ideal. 
The first theorem above can be reformulated as:
\begin{thm}\label{twideal}
Let $\mac J$ be the ideal of all planar graphs, and let $\mac I$ be an ideal of graphs. Then $\mac J\not\subseteq \mac I$ if and only if there exists $\ell$ such that every graph in $\mac I$ has tree-width at most $\ell$.
\end{thm}
This formulation makes its ``if-and-only-if'' nature more clear, and similar reformulations in terms of ideals exist for all rough structure theorems. 

It is helpful in this paper to designate one vertex of the tree $T$ as a {\em root}; and then a tree-decomposition $(T,(W_t:t\in V(T)))$ is called a {\em rooted tree-decomposition}. 
There are many more rough structure theorems; for instance, here are two that are well known (the {\em height} of a rooted tree-decomposition $(T,(W_t:t\in V(T)))$ is the length of the longest path in $T$ starting from the root):
\begin{thm}\label{path}
For every path $H$ there exists $m$ such that every graph not containing $H$ as a minor 
admits a rooted tree-decomposition of width at most $m$ and height at most $m$.
\end{thm}
\begin{thm}\label{erdos}
(Erd\H{o}s and P\'osa~\cite{erdos}:) If every component of $H$ has one vertex and one edge (a loop), there exists $\ell$ such that every graph not containing $H$ as a minor can be obtained 
from a forest by adding at most $\ell$ vertices (joined arbitrarily).
\end{thm}
And there are many more, some of which are substantially more complicated: for instance the ``graph minors structure theorem'' of Robertson and Seymour~\cite{GM16} is such a result, which we do  not state here.

In this series of papers we will look for rough structure theorems that lie between \ref{GM5} and \ref{GM1}; we want to exclude types of planar graph that are more complicated than forests, and see what we get.
A natural question is, is there a rough structure theorem for excluding a series-parallel graph? We have solved this, and will give the answer in a later paper, but it requires a good many preliminaries;
and there is a toy version of this question that we answer first, because the answer is much simpler and very pretty. The main goal of this paper is to prove a rough structure theorem that describes the 
series-parallel graphs that do not contain a given series-parallel graph as a minor.

Some terminology: if  $\mac T=(T,(W_t:t\in V(T)))$ is a tree-decomposition of $G$ and $t\in V(T)$, the {\em torso} of $\mac T$ at $t$ is the graph obtained from $G[W_t]$ by adding edges to make $W_s\cap W_t$
a clique for all $s\in V(T)$ adjacent to $t$. (More precisely, add a minimal set of edges to make $W_s\cap W_t$ a clique for all such $s$; do not introduce parallel edges.) If $e=st$ is an edge of $T$, the {\em adhesion} at $e$ is the set $W_s\cap W_t$. A tree-decomposition has {\em adhesion $k$} if all its adhesions have size exactly $k$.
(Many of the results in this paper are about tree-decompositions with adhesion two.)

\section{Statement of results}

Let us state our results. A planar drawing of a graph is {\em outerplanar} if there is a region that is incident with each vertex, and a graph is {\em outerplanar} if it admits an outerplanar drawing.
Rik Sengupta~\cite{sengupta} proved a result related to our question. Let us say an {\em apex-forest} is a graph obtained from a forest by adding one vertex, joined arbitrarily. (Thus, apex-forests are series-parallel.)
Sengupta proved:
\begin{thm}\label{senguptathm}
For every apex-forest $H$,
there exists $m$ such that every simple two-connected series-parallel graph with no $H$ minor admits a rooted tree-decomposition with adhesion two,
of height at most $m$, and such that each torso is outerplanar.
\end{thm}
This is a rough structure theorem, because it has a converse as usual: for all $m$, there is an apex-forest $H$ such that every two-connected series-parallel graph that admits such a tree-decomposition
contains no $H$ minor.

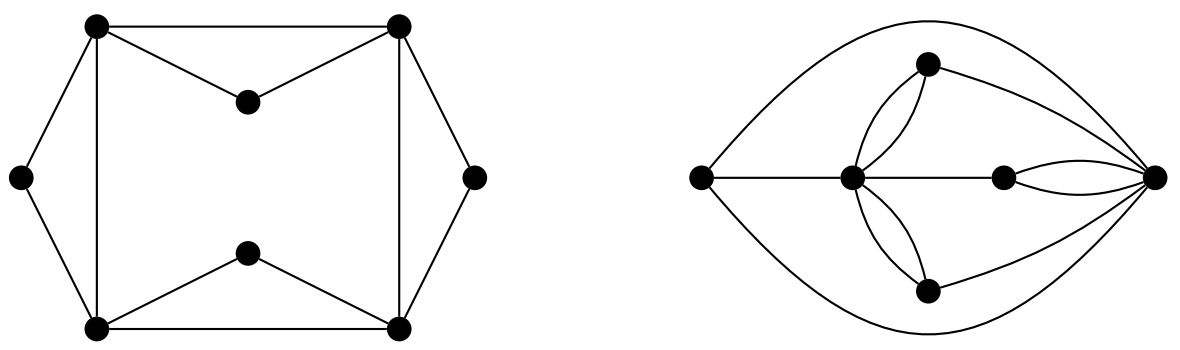
\begin{figure}[h]
\begin{center}
\begin{tikzpicture}[scale = 1]
\tikzstyle{every node}=[inner sep=1.5pt, fill=black,circle,draw]
\node (1) at (-1,1) {};
\node (2) at (1,1) {};
\node (3) at (1,-1) {};
\node (4) at (-1,-1) {};

\node (5) at (0,1/2) {};
\node (6) at (1.5,0) {};
\node (7) at (0,-1/2) {};
\node (8) at (-1.5,0) {};

\node (a) at (4,0) {};
\node (b) at (6,0) {};
\node (c) at (4.5,.75) {};
\node (d) at (5,0) {};
\node (e) at (4.5,-.75) {};
\node (f) at (3.0,0) {};

\foreach \to/\from in {1/2,2/3,3/4,1/4,1/5,2/5,2/6,3/6,3/7,4/7,4/8,1/8,a/d,a/f}
\draw[-] (\to)--(\from);
\foreach \to/\from in {b/c,e/b}
\draw[-] (\to) to [bend right = 10](\from);
\foreach \to/\from in {a/c,c/a,a/e,e/a}
\draw[-] (\to) to [bend right = 20](\from);

\foreach \to/\from in {d/b,b/d}
\draw[-] (\to) to [bend right = 20](\from);
\foreach \to/\from in {f/b,b/f}
\draw[-] (\to) to [bend right = 50, looseness = 1.5](\from);

\end{tikzpicture}
\end{center}
\caption{The left graph is outerplanar, but drawn with a non-outerplanar drawing. Its dual, on the right, is innerplanar but not an apex-forest.} \label{fig:innerplanar}
\end{figure}

If we take an outerplanar drawing of some connected (outerplanar) graph, and construct the dual graph, we obtain an apex-forest. But  a connected outerplanar graph in general admits planar drawings 
that are not outerplanar drawings, and, 
correspondingly, not all their duals are apex-forests. (See Figure \ref{fig:innerplanar}.) A connected graph is {\em innerplanar} if it admits a planar drawing
that is the dual of a drawing of some outerplanar graph. Let us say a disconnected graph is innerplanar if we can add edges to make a connected innerplanar graph.
(One can check that the innerplanar graphs are those that do not contain $K_4$ or $C_3^2$ as a minor, where $C_3^2$ is obtained from a cycle of length three by replacing each edge by two parallel edges.)

Our main result is:

\begin{thm}\label{SPthm}
For every series-parallel graph $H$,
there exists $m$ such that every simple two-connected series-parallel graph with no $H$ minor admits a rooted tree-decomposition with adhesion two, 
of height at most $m$, and such that each torso is either outerplanar or innerplanar.
\end{thm}

We will use the same machinery to show:
\begin{thm}\label{outerpthm}
For every outerplanar graph $H$,
there exists $m$ such that every simple two-connected series-parallel graph with no $H$ minor admits a rooted tree-decomposition with adhesion two, 
of height at most $m$, and such that each torso is innerplanar.
\end{thm}

\begin{thm}\label{innerpthm}
For every innerplanar graph $H$,
there exists $m$ such that every simple two-connected series-parallel graph with no $H$ minor admits a rooted tree-decomposition with adhesion two,
of height at most $m$, and such that each torso is either outerplanar or an apex-forest.
\end{thm}

\begin{figure}[h]
\begin{center}
\begin{tikzpicture}[scale = 1]
\tikzstyle{every node}=[inner sep=1.5pt, fill=black,circle,draw]
\def\h{.8}
\node (1) at ({cos(0)}, {sin(0)}) {};
\node (2) at ({cos(30)}, {sin(30)}) {};
\node (3) at ({cos(60)}, {sin(60)}) {};
\node (4) at ({cos(90)}, {sin(90)}) {};
\node (5) at ({cos(120)}, {sin(120)}) {};
\node (6) at ({cos(150)}, {sin(150)}) {};
\node (7) at ({cos(180)}, {sin(180)}) {};
\node (0) at (0,0) {};

\foreach \from/\to in {1/2,2/3,3/4,4/5,5/6,6/7,0/1,0/2,0/3,0/4,0/5,0/6,0/7}
\draw [-] (\from) -- (\to);
\end{tikzpicture}
\end{center}
\caption{A fan of length six.} \label{fig:fan}
\end{figure}
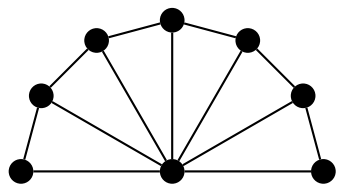

A {\em fan} is a simple graph obtained from a path by adding a vertex joined to each vertex of the path (see Figure \ref{fig:fan}), and its {\em length} is the length of the path. 
We will prove:
\begin{thm}\label{fanthm}
For every fan $H$, there exist $\ell, m$ such that every simple two-connected graph $G$ with no $H$ minor admits a
rooted tree-decomposition with adhesion two, of height at most $m$, and such that each torso is either a cycle or a three-connected graph with no path of length $\ell$.
\end{thm}
(Here, $G$ is a general graph, not necessarily series-parallel; but we can apply the same machinery because the result was already known for three-connected graphs, and all we need to do is to 
allow for 2-separations.)
It follows easily from \ref{path} (we omit the details) that
\begin{thm}\label{fanthm2}
For every fan $H$, there exist $\ell, m$ such that every simple two-connected graph $G$ with no $H$ minor admits a
rooted tree-decomposition of height at most $m$, and such that each torso either is a cycle or has at most $\ell$ vertices.
\end{thm}

The result \ref{SPthm} seemed very pretty, and that motivated us to ask, what if the graph we are decomposing has tree-width three rather than two? Does something similar happen?
It does indeed, although the proof is tougher.
We will show:
\begin{thm}\label{tw3}
For every series-parallel graph $H$, there exists $m$ such that if $G$ is a simple three-connected graph with tree-width at most three and with no $H$ minor, then 
$G$ admits a rooted tree-decomposition with adhesion three, of height at most $m$, and such that each torso is either:
\begin{itemize}
\item a near-outerplanar graph (defined later); or
\item obtained from an outerplanar graph by adding a vertex; or
\item obtained from an innerplanar graph by adding a vertex; or
\item a graph with path-width at most three.
\end{itemize}

\end{thm}

This can be extended to allow $G$ to be two-connected rather than three-connected, but we postpone that to a later paper since it requires some new methods.

Finally, we prove two other rough structure theorems, both easy, by proofs that do not use the method of ``scores'' that we used for all the earlier rough structure theorems. One is for excluding a ``stopsign'', and one is for excluding a ``caterpillar''. 

All these theorems have ``rough'' converses, and therefore are rough structure theorems. For instance, in the case of \ref{SPthm}, for all $m$ 
there is a series-parallel graph $H$ such that for every simple two-connected graph $G$, if $G$ admits 
a rooted tree-decomposition with adhesion two,
of height at most $m$, and such that each torso is either outerplanar or innerplanar, then $G$ has no $H$ minor. We sketch the proofs of the converses 
in Section~\ref{sec:converses}. 
(In most cases, given such a decomposition, one can delete the vertices of the bag of the root, and win by induction on $m$.)

\section{Innerplanar graphs}\label{sec:innerplanar}

There are several different ways to view innerplanar graphs, discussed in this section. 
Let $a,b$ be distinct vertices of a connected graph $C$, and let $a',b'$ be distinct vertices of a connected graph $C'$. Let $G$ be obtained from the disjoint union of $C,C'$ by identifying
$a=a'$ and $b=b'$, and let $H$ be obtained from the disjoint union of $C,C'$ by identifying
$a=b'$ and $b=a'$. The operation that takes $G$ to $H$ is called a {\em Whitney twist}. 

\begin{figure}[h]
\begin{center}
\begin{tikzpicture}[scale = 1]
\tikzstyle{every node}=[inner sep=1.5pt, fill=black,circle,draw]
\node (1) at (0,1.5) {};
\node (2) at (0,-1.5) {};
\node (3) at (-1.5,.5) {};
\node (4) at (-1.5,-.5) {};
\node (5) at (-.5,0) {};
\node (6) at (.5,0) {};
\node (7) at (1.5,.5) {};
\node (8) at (1.5,-.5) {};
\node (9) at (1.5,0) {};

\node (11) at (4,1.5) {};
\node (12) at (4,-1.5) {};
\node (13) at (2.5,.5) {};
\node (14) at (2.5,-.5) {};
\node (15) at (3.5,0) {};
\node (16) at (5.5,0) {};
\node (17) at (4.5,.5) {};
\node (18) at (4.5,-.5) {};
\node (19) at (4.5,0) {};

\foreach \to/\from in {1/3,1/5,1/6,1/7,2/4,2/5,2/6,2/8,3/4,3/5,4/5,6/7,6/8,7/9,8/9,6/9, 11/13,11/15,11/16,11/17,12/14,12/15,12/16,12/18,13/14,13/15,14/15,16/17,16/18,17/19,18/19,16/19}
\draw[-] (\to)--(\from);

\end{tikzpicture}
\end{center}
\caption{A Whitney flip.} \label{fig:flip}
\end{figure}
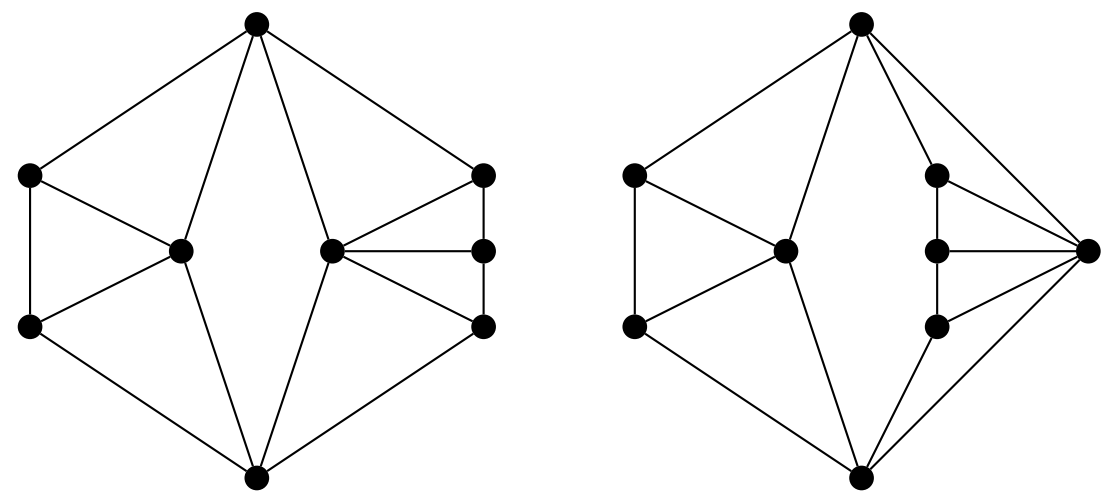

A {\em Whitney flip} is different; it is an operation that takes a planar drawing
of a graph to another drawing of the same graph, best shown by a figure. (See Figure \ref{fig:flip}.)

We start with a planar drawing (of a loopless two-connected graph, say, although the operation can be extended to other graphs), and choose a closed disc $\Delta$ in the plane whose boundary $F$ passes through two vertices $u,v$ and 
two regions of the drawing (and nothing else). 
Take a homeomorphism of $\Delta$ to itself that fixes $u,v$ and exchanges
the two components of $F\setminus \{u,v\}$, and replace the part of the drawing inside $\Delta$ by its image under this homeomorphism. This gives a new drawing of the same graph, and this operation
is called a {\em Whitney flip}. (An internet search shows that the two terms are often confused, and sometimes called a Whitney ``switch''.) Whitney proved that given 
two drawings of the same two-connected loopless graph in the sphere, we can take
one drawing to the other by a sequence of Whitney flips, up to a homeomorphism of the sphere to itself. (This is not explicitly stated in Whitney's papers as far as we know, but follows from results of~\cite{whitney}.)

\begin{thm}\label{flips}
A two-connected loopless graph is innerplanar if and only if it can be transformed by Whitney twists to an apex-forest.
\end{thm}
\Proof
The dual of a two-connected loopless drawing (not a cycle) is also two-connected and loopless (see \cite{whitney}). 
By Whitney's theorem, a drawing $D$ of a
two-connected loopless graph is a drawing of an outerplanar graph if and only if there is a sequence of drawings $D_1\LL D_n$ each obtained from its predecessor by a
Whitney
flip, with $D_1= D$ and $D_n$ an outerplanar drawing. Thus, a graph $G$ is innerplanar if and only if it has a drawing $D$ such that the dual drawing admits such a sequence; that is, 
(in terms of dual drawings),
if and only if there is a sequence of drawings 
$E_1\LL E_n$, with $E_1=D$ and $E_n$ the dual of an outerplanar drawing, where each $E_i$ is obtained from its predecessor by a Whitney twist (since 
a drawing is taken to a second drawing by a Whitney flip if and only if the dual of the first is obtained from the dual of the second by a Whitney twist). But the dual of a two-connected, loopless 
outerplanar drawing is a
two-connected loopless apex-forest. This proves \ref{flips}.~\bbox 

\begin{figure}[h]
\begin{center}
\begin{tikzpicture}[scale = 1]
\tikzstyle{every node}=[inner sep=1.5pt, fill=black,circle,draw]

\node (1) at (0,3) {};
\node (2) at (1,3) {};
\node (3) at (3,3) {};
\node (4) at (4,3) {};
\node (5) at (0,2) {};
\node (6) at (1,2) {};
\node (7) at (3,2) {};
\node (8) at (4,2) {};
\node (9) at (0,1) {};
\node (10) at (1,1) {};
\node (11) at (3,1) {};
\node (12) at (4,1) {};
\node (13) at (0,0) {};
\node (14) at (1,0) {};
\node (15) at (2,0) {};
\node (16) at (3,0) {};
\node (17) at (4,0) {};

\foreach \to/\from in
{1/6,2/5,3/8,4/7,5/10,6/9,6/11,7/10,7/12,8/11,9/14,10/13,11/15,11/17,12/16}
\draw[green, thick] (\to)--(\from);
\foreach \to/\from in
{1/5,2/6,3/7,4/8,5/6,7/8,6/10,7/11,9/10,10/11,11/12,9/13,10/14,10/15,11/16,12/17}
\draw[red, thick] (\to)--(\from);

\end{tikzpicture}
\end{center}
\caption{An innerplanar graph. The red tree is an inner tree.} \label{fig:innerpex}
\end{figure}
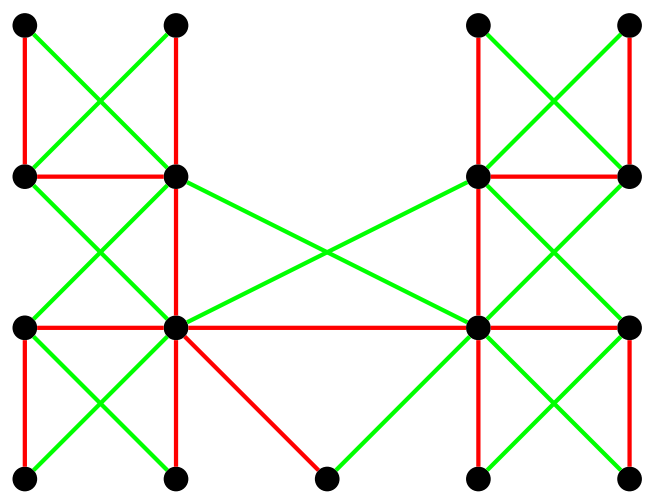

Let $G$ be a graph and $X\subseteq E(G)$; {\em adding a triangle} to $(G,X)$ results in $(G',X')$ where 
\begin{itemize}
\item $G'$ is obtained from $G$ by choosing an edge $uv$ of $G$ in $X$, and adding 
a new vertex $w$ adjacent to both $u,v$.  
\item $X'$ is obtained by choosing $f$ to be one of the
edges $uw,vw$, and setting $X'=X\cup \{f\}$.
\end{itemize}
(See Figure \ref{fig:innerpex} for a graph that can be grown in this way: start with the longest horizontal edge.)
\begin{thm}\label{growinn}
If there is a sequence $(G_1,X_1)\LL (G_n, X_n)$ where $G_1$ is a two-vertex graph with one edge (not a loop) and $X_1$ contains its edge, and $(G_i, X_i)$ is obtained from 
$(G_{i-1}, X_{i-1})$ by adding a triangle for each $i\ge 2$, then  $G_n$ is innerplanar.
\end{thm}
\Proof 
For $i\ge 2$, let $w_i$ be the vertex in $V(G_i)\setminus V(G_{i-1})$, and let $u_iv_i$ be the edge of $G_{i-1}$ that forms a triangle with $w_i$. Let $V(G_1) = \{w_0,w_1\}$. Then $V(G_n) = \{w_0\LL w_n\}$, and 
$$(\{w_0\LL w_{i-1}\},\{u_i, v_i, w_i\})$$
can be extended to a 2-separation of $G_n$, for each $i\ge 2$. By performing 
a Whitney twist at this 2-separation if necessary, for each $i\ge 2$, we can arrange that all edges in $X_n$ have a common end $w_0$; but then we have converted $G_n$ to an apex-forest by Whitney twists, and 
so $G_n$ is innerplanar, by \ref{flips}. This proves \ref{growinn}.~\bbox

Here is a third way to view innerplanar graphs. Take a tree $T$, and for each vertex $v\in V(T)$, take a tree $T_v$ with vertex set the set of neighbours of $v$ in $T$.
(The trees $T_v$ are necessarily edge-disjoint from $T$ and from each other.) Let $G$ be the union of $T$ and all the trees $T_v\;(v\in V(T))$. Then $G$ is innerplanar,
simple, and edge-maximal with these properties, and conversely, every edge-maximal simple innerplanar graph can be constructed in this way. If $G$ can be grown by adding 
edges to a tree $T$ in this way,  we call $T$ an {\em inner tree} for $G$.                                
(See Figure \ref{fig:innerpex}.)

\section{Rooted minors}

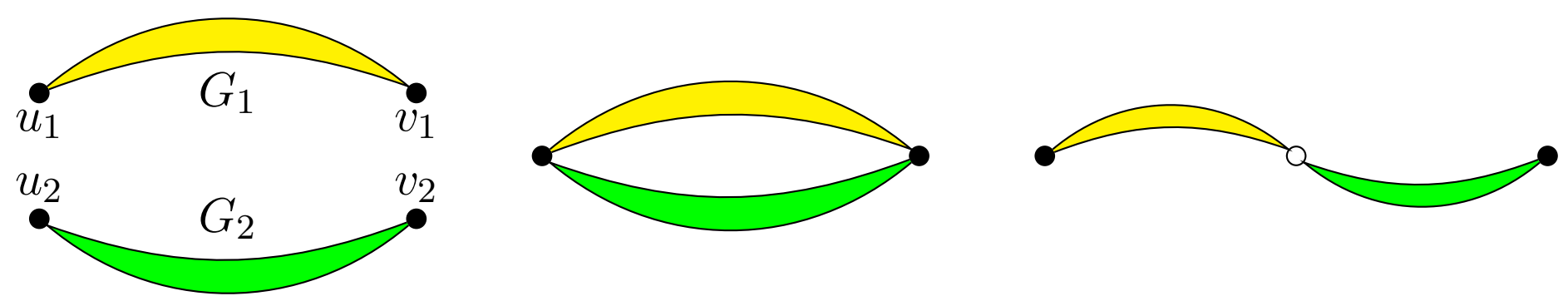
\begin{figure}[h]
\begin{center}
\begin{tikzpicture}[scale = 1]
\tikzstyle{every node}=[inner sep=1.5pt, fill=black,circle,draw]
\node (u) at (0,.5) {};
\node (v) at (3,.5) {};
\node (u2) at (0,-.5) {};
\node (v2) at (3,-.5) {};

\node (u') at (4,0) {};
\node (v') at (7,0) {};
\node (u'') at (8,0) {};
\node (w'') at (12,0) {};
\tikzstyle{every node}=[inner sep=1.5pt, circle,draw]
\node (v'') at (10,0) {};

\draw[fill=yellow] (u) to [bend left = 40] (v) to [bend right = 20] (u);
\draw[fill=green] (v2) to [bend left = 40] (u2) to [bend right = 20] (v2);

\draw[fill=yellow] (u') to [bend left = 40] (v') to [bend right = 20] (u');
\draw[fill=green] (v') to [bend left = 40] (u') to [bend right = 20] (v');

\draw[fill=yellow] (u'') to [bend left = 40] (v'') to [bend right = 20] (u'');
\draw[fill=green] (w'') to [bend left = 40] (v'') to [bend right = 20] (w'');

\tikzstyle{every node}=[]
\draw (u) node[below] {$u_1$};
\draw (v) node[below] {$v_1$};
\draw (u2) node[above] {$u_2$};
\draw (v2) node[above] {$v_2$};

\node at (1.5,.5) {$G_1$};
\node at (1.5,-.5) {$G_2$};

\end{tikzpicture}
\end{center}
\caption{Parallel and series connection.} \label{fig:doubling}
\end{figure}

A {\em $k$-rooted graph} $(G, r_1\LL r_k)$ consists of a graph $G$ and $k$ distinct vertices $r_1\LL r_k\in V(G)$. For the moment we will only need 
$k=2$. If $(G_1,u_1,v_1)$ and $(G_2,u_2,v_2)$ are 2-rooted graphs, their {\em parallel connection} is the 2-rooted graph $(G,u,v)$ where $G$ is obtained from the disjoint union of $G_1,G_2$ by 
identifying $u_1$ with $u_2$ (making $u$) and $v_1$ with $v_2$ (making $v$). (If $u_i,v_i$ are adjacent in $G_i$ for both $i = 1,2$, then these edges become 
parallel edges in the parallel connection.) The parallel connection of $(G,u,v)$ with itself is called its {\em parallel doubling}.
Similarly, the {\em series connection} of $(G_1,u_1,v_1)$ and $(G_2,u_2,v_2)$ is the 2-rooted graph $(G,u_1,v_2)$, where $G$ is obtained from the disjoint union of $G_1,G_2$ by identifying $v_1$ with $u_2$.
(The order matters: the series connection of $(G_2,u_2,v_2)$ and $(G_1,u_1,v_1)$ is different in general.)
The series connection of $(G,u,v)$ with itself is called its {\em series doubling}.
(See Figure \ref{fig:doubling}.)

For $k\ge 1$ we define a 2-rooted graph $\SP_k$ as follows. $\SP_1$ is the 2-rooted graph $(G,u,v)$ where $G$ consists of $u,v$ and an edge between them.
For even $k \ge 2$, $\SP_k$ is the series doubling of $\SP_{k-1}$, and for odd $k\ge 2$, $\SP_k$ is the parallel doubling of $\SP_{k-1}$. 
We show $\SP_1\LL \SP_5$ in Figure \ref{fig:SPk}.

\begin{figure}[h]
\begin{center}
\begin{tikzpicture}[scale = .9]
\tikzstyle{every node}=[inner sep=1.5pt, fill=black,circle,draw]
\def\w{-3}
\def\x{3}
\def\y{12}
\def\z{0}
\def\u{7}

\node (u) at (\w,0) {};
\node (v) at (\w+1,0) {};

\node (p) at (\z-1,0) {};
\node (r) at (\z+1,0) {};
\node (a) at (\x-1,0) {};
\node (d) at (\x+1,0) {};
\node (1a) at (\u-2,0) {};
\node (2a) at (\u+2,0) {};

\node (1) at (\y-2,0) {};
\node (2) at (\y+2,0) {};

\tikzstyle{every node}=[inner sep=1.5pt, circle,draw]
\node (b) at (\x,1) {};
\node (c) at (\x,.6) {};

\node (31) at (\u,1) {};
\node (51) at (\u-1.2,.95) {};
\node (61) at (\u-1,.65) {};

\node (71) at (\u+1.2,.95) {};
\node (81) at (\u+1,.65) {};

\node (q) at (\z,.6) {};
\node (3) at (\y,1) {};
\node (4) at (\y,2) {};
\node (5) at (\y-1.2,.95) {};
\node (6) at (\y-1,.65) {};
\node (7) at (\y+1.2,.95) {};
\node (8) at (\y+1,.65) {};
\node (9) at (\y-1.4,1.95) {};
\node (10) at (\y-1.2,1.65) {};
\node (11) at (\y+1.4,1.95) {};
\node (12) at (\y+1.2,1.65) {};
\foreach \to/\from in {u/v,p/q,q/r,a/b,a/c,d/b,d/c,1/5,1a/51,1/6,1a/61,3/5,31/51,3/6,31/61,3/7,31/71,3/8,31/81,2/7,2a/71,2/8,2a/81,1/9,1/10,4/9,4/10,4/11,4/12,2/11,2/12}
\draw[-] (\to)--(\from);

\end{tikzpicture}
\end{center}
\caption{$\SP_1\LL \SP_5$.} \label{fig:SPk}
\end{figure}

Thus, $\SP_k$ is a series-parallel graph (more exactly, if $(G,u,v)=\SP_k$ then $G$ is a series-parallel graph; we will usually not bother with this 
distinction); and it is easy to see that for every series-parallel graph $H$ there exists $k$ such that $H$ is a minor of $\SP_k$. 

If $G,H$ are graphs and $H$ is simple, a {\em model} of $H$ in $G$ is a function $\eta$ with domain $V(H)\cup E(H)$, such that
\begin{itemize}
\item for each $v\in V(H)$, $\eta(v)$ is a non-null connected subgraph of $G$, and $\eta(u),\eta(v)$ are vertex-disjoint for all distinct $u,v\in V(H)$;
\item for each edge $e\in E(H)$ with ends $u,v$, $\eta(e)$ is an edge of $G$ with one end in $V(\eta(u))$ and one end in $V(\eta(v))$. 
\end{itemize}
It is easy to see that $H$ is isomorphic to a minor of $G$ if and only 
if there is a model of $H$ in $G$.
If $(G,r_1\LL r_k)$ and $(H,s_1\LL s_k)$ are $k$-rooted graphs, we say the second is a {\em rooted minor} of the first if there is a model $\eta$
of $H$ in $G$ such that $r_i\in V(\eta(s_i))$ for $1\le i\le k$, and in this case we say that $\eta$ is a model of $(H,s_1\LL s_k)$ in $(G,r_1\LL r_k)$.

We will need to grow models of graphs into models of larger graphs, and the notation can be cumbersome. Here is a slight simplification: let $\eta$ be a model of $H$ in $G$, and let $P$ be a path of $G$
with one end in $\eta(u)$ for some $u\in V(H)$, and otherwise disjoint from each subgraph $\eta(v) (v\in V(H))$. Define $\eta'(u) = \eta(u)\cup P$, and $\eta'(v) = \eta(v)$ for all
$v\in V(H)\setminus \{u\}$, and $\eta'(e)=\eta(e)$ for each $e\in E(H)$. Then $\eta'$ is also a model of $H$ in $G$, and for brevity we say such a model is obtained from $\eta$ by ``contracting $P$''. We use similar terminology
in some more complicated situations, and we hope that the meaning will be clear.

\section{Refined and reduced tree-decompositions}

We need some more terminology and notation to handle a tree-decomposition $(T,(W_t:t\in V(T)))$. 
We 
direct all edges of $T$ away from the root. Thus, if we speak of an edge $t_1t_2$ of $T$, we imply that $t_1$ lies in the path of $T$ between the root and $t_2$. 
We call $t_1$
the
{\em parent} of $t_2$ (in $T$), and $t_2$ a {\em child} of $t_1$. More generally, if $s$ belongs to the path of $T$ between the root and $t$, we say $s$ is an {\em ancestor of $t$} and 
$t$ is a {\em descendant} of $s$.

A {\em triangle} means a cycle of length three. A graph is {\em $k$-connected} if it has at least $k+1$ vertices and there is no $X\subseteq V(G)$ with $|X|\le k-1$ such that $G\setminus X$ is not connected.
A rooted tree-decomposition is {\em reduced} if it satisfies:
\begin{itemize}
\item if $s,t\in V(T)$ are adjacent then $W_s\not\subseteq W_t$; and
\item for every path $P$ of $T$ with one end the root of $T$, no two edges of $P$ have the same adhesion. (To clarify, we are just requiring
that the two adhesions are not the same set; they are allowed to have the same size.)
\end{itemize}
It is {\em $k$-refined}, where $k\ge 1$ is an integer, if it has adhesion $k$, and each torso is either 
$(k+1)$-connected or a $(k+1)$-vertex complete graph.

We need:
\begin{thm}\label{goodtd}
Let $k\ge 1$, and let $G$ be a simple $k$-connected graph. Then $G$ admits a reduced $k$-refined rooted tree-decomposition.
\end{thm}
\Proof
Choose a tree-decomposition $(T,(W_t:t\in V(T)))$ with adhesion $k$, such that the set of all
pairs $\{u,v\}$ such that $\{u,v\}\subseteq W_t$ for some $t\in V(T)$ is minimal. We call this set its {\em set of nonseparated pairs}. Such a tree-decomposition is {\em acceptable.}
Choose a root $r$ for $T$. The {\em height-sum} of $T$ is the
sum of the lengths of all paths of $T$ that have the root $r$ as one end. Let $\mac T=(T,(W_t:t\in V(T)))$ be an acceptable rooted tree-decomposition such that the height-sum of $T$ is minimum.
We will show that $\mac T$ satisfies the theorem.

If there are adjacent
$s,t\in V(T)$ with
$W_s\subseteq W_t$, contract the edge $st$ to a vertex $u$ and define $W_u=W_t$. This operation does not change the set of nonseparated pairs; and so it gives another acceptable tree-decomposition
with smaller height-sum, a contradiction.
Thus $\mac T$ satisfies the first bullet in the definition of ``reduced''.

If there are distinct edges $e_1,e_2$ of $T$ with the same adhesion, and $e_1$ is
between $e_2$ and the root $r$, let $e_i = s_it_i$ for $i = 1,2$, where there is a path of $T$ passing through $r,s_1,t_1,s_2,t_2$ in order
(possibly $r=s_1$ or $s_2 = t_1$). Let $T'$ be the tree obtained from $T$ by deleting $e_2$ and adding a new edge $s_1t_2$. Then
$(T', (W_t:t\in V(T)))$ is another acceptable tree-decomposition with smaller height-sum, again a contradiction. So $\mac T$ satisfies the second bullet,
and hence is reduced.

To show that $\mac T$ is $k$-refined, let
$t\in V(T)$, and suppose that the torso of $\mac T$
at $t$ is not $(k+1)$-connected and not a $(k+1)$-vertex complete graph. If this torso is a complete graph, then $|W_t|\le k$, and since all adhesions have size 
$k$ and $\mac T$ is reduced, it follows that $t$ has degree zero in $T$, and so $V(G)=W_t$, a contradiction since $G$ is $k$-connected. Thus, the torso at $t$
is not complete, and consequently there exists $X\subseteq W_t$ with $|X|\le k$, and a partition
$(Y_1,Y_2)$ of $W_t\setminus X$ such that $Y_1,Y_2\ne \emptyset$ and there are no edges of the torso between $Y_1,Y_2$. Hence we may partition
the set of edges of $T$ incident with $t$ into two sets $E_1,E_2$, such that for $i = 1,2$, if $st\in E_i$ then
$W_s\cap W_t\subseteq Y_i\cup X$. Since $G$ is $k$-connected, it follows that $|X|=k$.  Let $T'$ be obtained from $T$ by: deleting $t$; adding two new vertices $t_1,t_2$, adjacent to each other;
and for $i = 1,2$ and each $st\in E_i$, adding an edge between $s,t_i$. Then $T'$ is a tree; define $W_{t_i}=X\cup Y_i$ for $i = 1,2$.
It follows that $(T',(W_t:t\in V(T')))$ is a tree-decomposition of $G$, all its adhesions have size $k$, and its set of
nonseparated pairs is strictly reduced, a contradiction. This proves that $\mac T$ is $k$-refined, and so 
proves \ref{goodtd}.~\bbox

If $(T,(W_t:t\in V(T)))$ is a rooted tree-decomposition of $G$, and 
$e=t_1t_2\in E(T)$, then $T\setminus e$ has two components $T_1,T_2$ say, where the root belongs to $V(T_1)$ (and hence $t_i\in V(T_i)$ for $i = 1,2$).
We define $B_e=\bigcup_{t\in V(T_1)}W_t$ and $A_e=\bigcup_{t\in V(T_2)}W_t$ (``before'' $e$ and ``after'' $e$); thus, $(A_e,B_e)$ is a separation of $G$, and $A_e\cap B_e$ is the adhesion at $e$.

If $T$ is a tree, a {\em leaf} of $T$ is a vertex with degree one in $T$, and a {\em leaf-edge} is an edge incident with a leaf.
If $T$ is a tree with root $r$, and $S$ is a subtree of $T$, the {\em root} of $S$ is defined to be the vertex of $S$ closest to $r$ in $T$. If the root of $S$ is a leaf of $S$, the corresponding leaf-edge is 
a {\em root-edge} for $S$. A vertex of $S$ is {\em internal} if it has degree at least two in $S$. If $S$ has an internal vertex, $S^-$ is the subtree induced on the set of all internal vertices of $S$.
Let $(T,(W_t:t\in V(T)))$ be a rooted tree-decomposition of some graph $G$; if $S$ is a subtree of $T$ with a root-edge and an internal vertex, we define
the {\em $S$-torso} $H(S)$ to be the graph defined as follows. Let $G_1$ be the union of $G[W_s]$ over 
all $s\in V(S^-)$. The 
$S$-torso is obtained from $G_1$ by adding new edges joining each two members of $A_e\cap B_e$ (if they are currently nonadjacent), for each edge $e\in E(T)$ that has exactly one end in $V(S^-)$.

Here is a general proof method that we will use repeatedly. 
We have some $k$-connected graph $G$, and a rooted $k$-refined tree-decomposition $(T,(W_t:t\in V(T)))$ for $G$.
Each edge $e\in E(T)$
has a ``score'' $\sc(e)$. Scores are ordered by some partial order (in this paper, scores will be integers ordered in their natural order, but the method works more generally). For each path $P$ of $T$ with one end 
the root of $T$, the scores of the edges of $P$ monotonically (non-strictly) decrease with distance from the root, and take only a bounded number of different values, say at most $\ell$ values. 
Let $S$ be a subtree of $T$ such that:
\begin{itemize}
\item $S$ has a root-edge;
\item every edge of $S$ has the same score; and
\item $S$ is maximal with these properties.
\end{itemize}
We call such a subtree $S$ a {\em plateau}. Every edge of $T$ therefore belongs to a
unique plateau (because if two plateaus share an edge, their union also satisfies the first two bullets above, contrary to their maximality).
It follows that every path of $T$ with one end the root contains edges of at most $\ell$ plateaus.

Let $F$ be the set of all edges of $T$ that are leaf-edges of some plateau (including the root-edge of each plateau). Thus, each component of $T\setminus F$ is either a singleton vertex,
or obtained from some plateau by deleting its leaves. Let $T'$ be the rooted tree obtained from $T$ by contracting each component $C$ of $T\setminus F$ to a vertex $t_C$, and
setting the root of $T'$ to be the root of $T$ (which is a singleton component of $T\setminus F$). Thus, $E(T') = F$. For each component
$C$ of $T\setminus F$:
\begin{itemize}
\item if $|V(C)|=1$, let $V(C)=\{t\}$ say, and define $W'_{t_C} = W_t$;
\item if $|V(C)|>1$, let $C$ be obtained from a plateau $S$ by deleting its leaves, and define $W'_{t_C}=\bigcup_{t\in V(S^-)}W_t$.
\end{itemize}
Then $(T',(W'_t:t\in V(T')))$ is a rooted tree-decomposition of $G$ with adhesion $k$, and its height is at most $2\ell$, and each of its torsos is either a torso of some vertex of $T$ in 
$(T,(W_t:t\in V(T)))$ (and hence either $(k+1)$-connected or a $(k+1)$-vertex complete graph), or equal to $H(S)$ for some plateau $S$
with an internal vertex.
We have proved:
\begin{thm}\label{genthm}
Let $k,\ell\ge 1$, let $G$ be a $k$-connected graph, and let $(T,(W_t:t\in V(T)))$ be a $k$-refined tree-decomposition for $G$. For each edge $e\in E(T)$, let $\sc(e)$ be an integer,
the ``score'' of $e$.
Suppose that for every path $P$ of $T$ with one end the root of $T$:
\begin{itemize}
\item if $e,e'\in E(P)$ and $e$ is between $e'$ and the root, then $\sc(e')\le \sc(e)$; and
\item the edges of $P$ have at most $\ell$ different scores.
\end{itemize}
Then $G$ admits a rooted tree-decomposition with adhesion $k$ and of height at most $2\ell$, and each of its torsos is either $(k+1)$-connected, or a $(k+1)$-vertex complete graph,
or the $S$-torso $H(S)$ of some plateau $S$ with an internal vertex.
\end{thm}

\section{The proof of \ref{SPthm}}

In this section we prove \ref{SPthm}, which we restate:
\begin{thm}\label{SPthmagain}
For every series-parallel graph $H$,
there exists $m$ such that every simple two-connected series-parallel graph with no $H$ minor admits a rooted tree-decomposition with adhesion two,
of height at most $m$, and such that each torso is either outerplanar or innerplanar.
\end{thm}
\Proof
We may assume that $H=\SP_k$ for some choice of $k$. Let $G$ be a simple two-connected series-parallel graph with no $\SP_k$ minor. Choose a 
$2$-refined reduced rooted tree-decomposition
as in \ref{goodtd}. Since $G$ is series-parallel, no torso is three-connected; and so $|W_t|=3$ for each $t\in V(T)$.
For each $e\in E(T)$, we say the {\em score} $\sc(e)$ of $e$ is the maximum $d$ such that $\SP_d$ is a rooted minor of $(G[A_e],u,v)$, where
$A_e\cap B_e = \{u,v\}$. (We will have to change the definition of ``score'' to prove \ref{outerpthm} and \ref{innerpthm}, but this is the version we need for \ref{SPthm}.)
From the symmetry of $\SP_d$ exchanging its roots, the order of $u,v$ makes no difference.

If $e,f\in E(T)$ and $e$ belongs to the path of $T$ between the root and $f$, then the score of $e$ is at least the score of $f$, because there are two disjoint paths between 
$A_e\cap B_e$ and $A_f\cap B_f$, so if $\SP_d$ is a rooted minor of $(G[A_f],u,v)$ (where $A_f\cap B_f = \{u,v\}$) then it is also a rooted minor of $(G[A_e],w,x)$ where 
$A_e\cap B_e=\{w,x\}$, as we see by contracting the two paths. 
Thus, for every path of $T$ with one end the root, the scores of its edges (non-strictly) decrease monotonically with distance from the root.
Since every edge of $T$ has score between $1$ and $k-1$ (because $\SP_k$ is not a minor of $G$), it follows that every path of $T$ with one end the root contains edges of at most $k$ plateaus. 
By \ref{genthm}, 
$G$ admits a rooted tree-decomposition with adhesion two and of height at most $2k$, and each of its torsos is either a triangle,
or the $S$-torso $H(S)$ of some plateau $S$ with an internal vertex.
Hence, to complete the proof of \ref{SPthm}, it remains only to check that each $S$-torso $H(S)$ is outerplanar or innerplanar. 

Thus, let $S$ be a plateau with an internal vertex and with root $r_0$, and let $d$ be the common score of the edges of $S$.
Let $s\in V(S)\setminus \{r_0\}$, and let $r$ be its parent (thus $r\in V(S)$) and $f=rs$. Let $W_s=\{v_1,v_2,v_3\}$, where $A_f\cap B_f = \{v_1,v_2\}$. Since $(T,(W_t:t\in V(T)))$ is reduced, every edge 
$e\ne f$
of $T$ incident with $s$ satisfies $A_e\cap B_e=\{v_1,v_3\}$ or  $A_e\cap B_e=\{v_2,v_3\}$.
\\
\\
(1) {\em If $d$ is even, there is at most one edge $e\in E(S)$ incident with $s$ with $A_e\cap B_e = \{v_1,v_3\}$, and the same for $\{v_2,v_3\}$.}
\\
\\
Since $d$ is even, $\SP_{d+1}$ is the parallel doubling of $\SP_d$. 
Suppose that $e,e'\in E(T)$ are distinct edges of $S$ both incident with $s$ and with 
$A_e\cap B_e = A_{e'}\cap B_{e'} = \{v_1,v_3\}$. Since $e\in E(S)$, and so $e$ has score $d$, it follows that $\SP_d$ is a rooted minor of $(G[A_e],v_1,v_3)$, and similarly 
it is a rooted minor of $(G[A_{e'}],v_1,v_3)$. By combining the corresponding two models, it follows that the parallel doubling $\SP_{d+1}$ of $\SP_d$ is a rooted minor of $(G[A_e\cup A_{e'}],v_1,v_3)$.
Let $\eta$ be the corresponding model of $\SP_{d+1}$.
Since $G$ is two-connected, there is a path of $G\setminus \{v_1\}$ between $v_2,v_3$, say $P$, and by contracting this path we convert $\eta$
to a model of $\SP_{d+1}$ in $G[A_f]$; and so $\SP_{d+1}$ is a rooted minor of $(G[A_f],v_1,v_2)$, a contradiction.
This proves (1).
\\
\\
(2) {\em If $d$ is odd, there do not exist $e,e'\in E(S)$ incident with $s$ with $A_e\cap B_e = \{v_1,v_3\}$ and  $A_{e'}\cap B_{e'} = \{v_2,v_3\}$.}
\\
\\
Since $d$ is odd, $\SP_{d+1}$ is the series doubling of $\SP_d$. If such $e,e'$ exist, then by combining the corresponding two models, we obtain a model of $\SP_{d+1}$ in $(G[A_f],v_1,v_2)$, a contradiction.
This proves (2).

\bigskip

Every edge of the $S$-torso has both ends in $W_s$ for some $s\in V(S)\setminus \{r_0\}$. Let $J$ be the graph with vertex set $\bigcup_{s\in V(S)\setminus \{r_0\}}W_s$ and with edges every pair
of distinct vertices that both belong to some $W_s$ for some $s\in V(S)\setminus \{r_0\}$. If $d$ is even then from (1), it follows that 
$J$ is outerplanar, and hence so is the $S$-torso.  If $d$ is odd then from (2) and \ref{growinn}, it follows that $J$ is innerplanar, and hence so is the $S$-torso. 
This proves \ref{SPthmagain}.~\bbox

\section{When the $S$-torso is a cycle}

To prove \ref{outerpthm} and \ref{fanthm}, we need the following lemma:
\begin{thm}\label{getcycle}
Let $G$ be a two-connected graph, and let $(T,(W_t:t\in V(T)))$ be a reduced 2-refined rooted tree-decomposition of $G$. Let $S$ be a subtree of $T$ with a root-edge $e_0=r_0s_0$ and an internal vertex, and let
$H(S)$ denote the $S$-torso. If $H(S)$ is not a cycle, there are edges $e,e'$ of $S$, such that
\begin{itemize}
\item the head of $e$ is the tail of $e'$;
\item there are two disjoint paths $P_1,P_2$
of $G$ between $A_{e}\cap B_{e}$ and $A_{e'}\cap B_{e'}$; and
\item there is a path of $G[A_{e}\cap B_{e'}]$ between $V(P_1), V(P_2)$.
\end{itemize}
\end{thm}
\Proof Suppose there are no such $e,e'$ (we call this ``the supposition''). We will prove that $H(S)$ is a cycle, by induction on $|V(S)|$.
\\
\\
(1) {\em $|W_s|=3$ for each $s\in V(S^-)$.}
\\
\\
Let $s\in V(S^-)$. Thus $s\ne r_0$, and is not a leaf of $S$. Choose edges $e,e'$ of $S$ incident with $s$, such that $e$ has head $s$ and $e'$ has tail $s$. Let $H$ be the
torso of $(T,(W_t:t\in V(T)))$ at $s$.
Let $X=A_e\cap B_e$ and $X'=A_{e'}\cap B_{e'}$; then there are edges $f_e,f_{e'}$ of $H$ joining the two vertices in $X$ and joining the two vertices in $X'$ respectively.
Suppose that $H$ is three-connected. Then $|V(H)|\ge 4$, so deleting two edges of $H$ does not disconnect $H$; and hence $H\setminus \{f_e,f_{e'}\}$ is connected.
Since $G$ is two-connected, there are two disjoint paths $Q_1,Q_2$
of $H$ between $X,X'$. Since $H\setminus \{f_e,f_{e'}\}$ is connected, there is a path $Q_0$ of $H\setminus \{f_e,f_{e'}\}$ between $V(Q_1),V(Q_2)$ with no internal vertex
in $V(Q_1)\cup V(Q_2)$. All vertices of these paths belong to $W_s$, but some of their
edges might not be edges of $G$, because $H$ is not necessarily a subgraph of $G$. 
For each edge $g$ of $Q_0, Q_1$ or $Q_2$
that is not an edge of $G$,
there is an edge $h$ of $T$ incident with $s$ with $A_h\cap B_h$ equal to the set of ends of $g$, and there is a path $R_g$ of $G[A_h]$ joining the two ends of $g$.
By replacing each edge $g$ of $Q_0$ that is not in $E(G)$ by the path $R_g$, we obtain a path $P_0$ of $G$ joining the ends of $Q_0$. Define $P_1,P_2$ similarly; then $P_0,P_1, P_2$
violate the supposition. Thus $H$ is not three-connected; and so $|W_s|=3$ since $(T,(W_t:t\in V(T)))$ is 2-refined. This proves (1).

\bigskip

We recall that $r_0s_0$ is the root-edge of $S$.
Since $S$ has an internal vertex, it follows that $s_0$ is internal, and so $|W_{s_0}|=3$ by (1).  Let $W_{s_0}= \{v_1,v_2,w\}$, where
$\{v_1,v_2\}$ is the adhesion at $e_0$. For each edge $e\ne e_0$ of $T$ incident with $s_0$, the adhesion at $e$
is one of the sets $\{v_1,w\},\{v_2,w\}$, since no such edge has adhesion $\{v_1,v_2\}$ by the second property of ``reduced''.
Let $E_i$ be the set of edges of $T$ incident with $s_0$ with adhesion
$\{v_i,w\}$, for $i = 1,2$.
For each edge $e\in E_1\cup E_2$, since $G$ is two-connected and $A_e\ne A_e\cap B_e$ (by the first property of ``reduced''), there is a path of
$G[A_e]$ between the two members of $A_e\cap B_e$, say $P_e$.

Suppose that there is an edge $e_1$ of $S$ that belongs to $E_i$, for some $i\in \{1,2\}$. By the supposition, applied to the edges $e_0,e_1$,
there is no path in
$G[A_{e_0}\cap B_{e_1}]$ between $v_i,w$. Consequently, $v_i,w$ are nonadjacent in $G$; and moreover,
there is no path $P_e$ for $e\in E_i\setminus \{e_1\}$ and so $E_i=\{e_1\}$. In summary, either $E_i\cap E(S)=\emptyset$, or $|E_i|=1$ and $E_i\subseteq E(S)$.
In either case, at most one edge of $S$ belongs to $E_i$, so $s_0$ has degree at most three in $S$.

Since $|E(S)|\ge 2$, at least one of $E_1,E_2$ contains an edge of $S$; so we may assume that  $|E_1|=1$ and $E_1\subseteq E(S)$, say.
For $i = 1,2$, let $F_i$ be the set of $e\in E_i\cap E(S)$ such that both ends of $e$ are internal vertices of $S$.
If $F_1=F_2=\emptyset$, then the $S$-torso $H(S)$ has vertex set $W_{s_0}$; and since $|E_1|=1$, the edge $v_1w$ of $H(S)$ is not parallel to any other edge of $H(S)$, and therefore $H(S)$ is a cycle.
Thus we may assume that $F_1\ne \emptyset$. Now there are three cases, depending on whether no edge of $S$ belongs to $E_2$, or some edge of $S$ is in $E_2$ but $F_2=\emptyset$, or $F_2\ne \emptyset$.

If no edge of $S$ belongs to $E_2$, then $s_0$ has degree two in $S$; let $S_1$ be the tree $S\setminus \{r_0\}$. If some edge of $S$ belongs to $E_2$, let $S_1,S_2$
be the maximal subtrees of $S\setminus \{r_0\}$ with leaf-edges $e_1,e_2$, where $e_i$ is the unique edge of $S$ in $E_i$ for $i = 1,2$. From the inductive hypothesis, applied to $S_1$, the
$S_1$-torso $H(S_1)$ is a cycle, and it contains the edge $f_{e_1}$ joining the two vertices of $A_{e_1}\cap B_{e_1}$. But $v_1,w$ are nonadjacent in $G$, and
$E_1=\{e_1\}$, so the edge $f_{e_1}$ does not belong to $H(S)$. Hence, if $F_2=\emptyset$, then $H(S)$ is obtained from the cycle $H(S_1)$ by deleting $f_{e_1}$ and adding $v_2$ and two new edges
$v_1v_2,v_2w$; and so it is a cycle. Thus we may assume that $F_2\ne\emptyset$. But then both $H(S_1), H(S_2)$ are cycles, and the edges $f_{e_1},f_{e_2}$ do not belong to $H(S)$, and so
$H(S)$ is obtained from the union of $H(S_1),H(S_2)$ by deleting $f_{e_1},f_{e_2}$ and adding $v_1v_2$. This proves \ref{getcycle}.~\bbox

\section{The proof of \ref{outerpthm}}

Next we prove \ref{outerpthm}. The proof is similar to that of \ref{SPthm}, except we have to change the definition of ``score'' and use \ref{getcycle}. Let us define
$\OP_1=\SP_1$, and for $d\ge 2$ even, let $\OP_d$ be the series doubling of $\OP_{d-1}$.
For $d\ge 2$ odd, let $\OP_d$ be obtained from $\OP_{d-1}$ by adding an edge between its two roots.  (See Figure \ref{fig:OPk}.)

\begin{figure}[h]
\begin{center}
\begin{tikzpicture}[scale = .9]
\tikzstyle{every node}=[inner sep=1.5pt, fill=black,circle,draw]
\def\w{-3}
\def\x{3}
\def\y{12}
\def\z{0}
\def\u{7}

\node (u) at (\w,0) {};
\node (v) at (\w+1,0) {};

\node (p) at (\z-1,0) {};
\node (r) at (\z+1,0) {};

\node (p1) at (\x-1,0) {};
\node (r1) at (\x+1,0) {};

\node (a) at (\u-2,0) {};
\node (e) at (\u+2,0) {};

\node (a1) at (\y-2,0) {};
\node (e1) at (\y+2,0) {};

\tikzstyle{every node}=[inner sep=1.5pt, circle,draw]
\node (q) at (\z,.6) {};
\node (q1) at (\x,.6) {};

\node (b) at (\u-1,1) {};
\node (c) at (\u,.5) {};
\node (d) at (\u+1,1) {};
\node (b1) at (\y-1,1) {};
\node (c1) at (\y,.5) {};
\node (d1) at (\y+1,1) {};

\foreach \to/\from in {u/v,p/q,q/r,p1/q1,q1/r1,p1/r1, a/b,a/c,b/c,d/c,d/e,c/e, a1/b1,a1/c1,b1/c1,c1/d1,c1/e1,d1/e1,a1/e1}
\draw[-] (\to)--(\from);

\end{tikzpicture}
\end{center}
\caption{$\OP_1\LL \OP_5$.} \label{fig:OPk}
\end{figure}
We restate \ref{outerpthm}:
\begin{thm}\label{outerpthmagain}
For every outerplanar graph $H$,
there exists $m$ such that every simple two-connected series-parallel graph with no $H$ minor admits a rooted tree-decomposition with adhesion two,  
of height at most $m$, and such that each torso is innerplanar.
\end{thm}
\Proof
Every outerplanar graph is a minor of $\OP_\ell$ for some $\ell$, so we may assume that $H =\OP_\ell$ for some $\ell$.
Let $G$ be a two-connected simple series-parallel graph that does not contain $\OP_\ell$ as a minor. 
Take a reduced 2-refined tree-decomposition $(T,(W_t:t\in V(T)))$ of $G$ of width two as before. 
Now the score of each edge $e$ of $T$ is defined as the maximum $d$ such that $\OP_d$ is a rooted minor of $(G[A_e],u,v)$,
where $A_e\cap B_e=\{u,v\}$. As before, the order of $u,v$ here does not matter. Define plateaus as before; so $G$ admits  a rooted tree-decomposition with adhesion two, 
with height at
most $2\ell$; and each of its torsos
is either a triangle, or the $S$-torso for some plateau $S$ with an internal vertex. Hence, to complete the proof of \ref{outerpthm}, 
it remains only to check that for each plateau $S$ with an internal vertex, the $S$-torso is 
innerplanar. Let $d$ be the common score of the edges of $S$. 

Suppose first that $d$ is even. 
If $e,e'$ are edges of $S$ and the head of $e$ equals the tail of $e'$, then since $e,e'$ have the same score $d$ (which is even), the conclusion of \ref{getcycle} is 
false. 
Hence, by applying \ref{getcycle} to $S$, we deduce that the $S$-torso is a cycle and hence innerplanar.

Now suppose that $d$ is odd. Then again the $S$-torso is innerplanar, by the same argument as in the proof of \ref{SPthmagain}. This proves \ref{outerpthmagain}.~\bbox
\section{The proof of \ref{fanthm}}

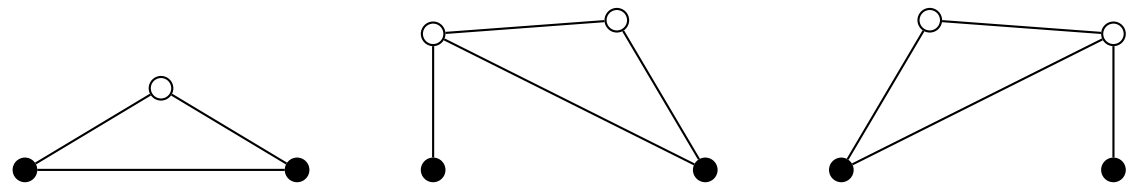
\begin{figure}[h]
\begin{center}
\begin{tikzpicture}[scale = .9]
\tikzstyle{every node}=[inner sep=1.5pt, fill=black,circle,draw]
\node (p) at (-1,0) {};
\node (r) at (+1,0) {};

\node (w1) at (2,0) {};
\node (r1) at (4,0) {};

\node (p2) at (5,0) {};
\node (w2) at (7,0) {};

\tikzstyle{every node}=[inner sep=1.5pt, circle,draw]
\node (q) at (0,.6) {};
\node (p1) at (2,1) {};
\node (q1) at (3.35,1.1) {};

\node (r2) at (7,1) {};
\node (q2) at (5.65,1.1) {};
\foreach \to/\from in {p/q,q/r,p/r, p1/q1,p1/r1,q1/r1,w1/p1,p2/q2,p2/r2,q2/r2,r2/w2}
\draw[-] (\to)--(\from);

\end{tikzpicture}
\end{center}
\caption{Leggings.} \label{fig:legs}
\end{figure}

If $(G,u,v)$ is a 2-rooted graph, its {\em left-legging} is the series connection
of $\SP_1$ and  $(G,u,v)$, that is, the rooted graph $(H,w,v)$, where $H$ is obtained from $G$ by adding a new vertex $w$ and the edge $uw$. Similarly,
its {\em right-legging} is the series connection of  $(G,u,v)$ and $\SP_1$. (See Figure \ref{fig:legs}.)

Let us define $\F_1 = \SP_1$; for $d\ge 2$, even, let $\F_d$ be the right-legging of $\F_{d-1}$, and for $d\ge 2$, odd, let $\F_d$ be obtained from $\F_{d-1}$ by adding an edge between
the roots of $\F_{d-1}$. (See Figure \ref{fig:fans}.)
\begin{figure}[h]
\begin{center}
\begin{tikzpicture}[scale = .9]
\tikzstyle{every node}=[inner sep=1.5pt, fill=black,circle,draw]
\node (u1) at (0,0) {};
\node (v1) at (1,0) {};
\node (u2) at (2,0) {};
\node (v2) at (3,0) {};
\node (u3) at (4,0) {};
\node (v3) at (5,0) {};

\node (u4) at (6,0) {};
\node (v4) at (7,0) {};
\node (u5) at (8,0) {};
\node (v5) at (9,0) {};
\tikzstyle{every node}=[inner sep=1.5pt, circle,draw]
\node (p2) at (3,.5) {};
\node (p3) at (5,.5) {};
\node (p4) at (7,.5) {};
\node (p5) at (9,.5) {};
\node (q4) at (7,1) {};
\node (q5) at (9,1) {};

\foreach \to/\from in {u1/v1,u2/p2,v2/p2,u3/p3,v3/p3,u3/v3, u4/p4,u4/q4,v4/p4,p4/q4,u5/v5,u5/p5,u5/q5,v5/p5,p5/q5}
\draw[-] (\to)--(\from);

\end{tikzpicture}
\end{center}
\caption{$\F_1\LL \F_5$.} \label{fig:fans}
\end{figure}

We also need the following theorem, proposition 3.8 of Ding, Dziobiak, and Wu~\cite{wheel}:
\begin{thm}\label{3connwheel}
For every $k$ there exists $\ell$ such that every three-connected graph with a path of length at least $\ell$ contains a $k$-wheel (and hence $\F_k$) as a minor.
\end{thm}

Now we prove \ref{fanthm}, which we restate:
\begin{thm}\label{fanthmagain}
For every fan $H$, there exist $\ell, m$ such that every simple two-connected graph $G$ with no $H$ minor admits a
rooted tree-decomposition with adhesion two, of height at most $m$, and such that each torso is either a cycle or a three-connected graph with no path of length $\ell$.
\end{thm}
\Proof
Let $G$ be a simple two-connected graph not containing $\F_k$ as a minor. Take a reduced 2-refined rooted tree-decomposition $(T, (W_t:t\in V(T)))$ as in \ref{goodtd}.
For each edge $e$ of $T$, let $A_e\cap B_e = \{u,v\}$; we define the {\em score} of $e$ to be $d_1+d_2$ where  $d_1$ is maximum  such that 
$\F_{d_1}$ is a rooted minor of $(G[A_e],u,v)$, and $d_2$ is maximum such that $\F_{d_2}$ is a rooted minor of $(G[A_e],v,u)$. Thus, all scores are less than $2k$. By \ref{genthm}, $G$ admits a rooted tree-decomposition 
with adhesion two, and of height at most $4k$, such that each of its torsos is three-connected (and therefore has no path of length $\ell$, where $\ell$ is as in \ref{3connwheel}), or a triangle,  or equal to $H(S)$ for some plateau $S$ with an internal vertex. 
Thus, it remains to show that for every plateau $S$ 
with an internal vertex, $H(S)$ is a cycle. Let $S$ be such a plateau, and let $d$ be the common score of the edges of $S$.  Let $e_0$ be the root-edge of $S$, let $A_{e_0}\cap B_{e_0} = \{u_0,v_0\}$, and 
choose $d_1,d_2$ as in the definition of score, with $d_1+d_2=d$. 
\\
\\
(1) {\em For every edge $e$ of $S$, let $P_1,P_2$ be disjoint paths between the adhesions of $e_0$ and of $e$, where $u_0\in V(P_1)$. Let the adhesion at $e$ be $\{u,v\}$, where $u\in V(P_1)$.
Then $d_1$ is maximum such that $\F_{d_1}$ is a rooted minor of $(G[A_e],u,v)$, and $d_2$ is maximum such that $\F_{d_2}$ is a rooted minor of $(G[A_e],v,u)$.}
\\
\\
Let $d_1',d_2'$ be maximum such that 
$\F_{d_1'}$ is a rooted minor of $(G[A_e],u,v)$, and $\F_{d_2'}$ is a rooted minor of $(G[A_e],v,u)$. Thus $d_1'+d_2'=d$, since $e,e_0$ have the same score. But
since $\F_{d_1'}$ is a rooted minor of $(G[A_e],u,v)$, by contracting the paths $P_1,P_2$ we deduce that $\F_{d_1'}$ is a rooted minor of $(G[A_{e_0}],u_0,v_0)$, and so $d_1'\le d_1$. 
Similarly $d_2'\le d_2$; but $d_1+d_2=d_1'+d_2'$, and the claim follows. This proves (1).

\bigskip

If one of $d_1,d_2$ is even, then it follows immediately from \ref{getcycle} that $H(S)$ is a cycle. If they are both odd, then since $S$ has an internal vertex, it has two edges $e,e'$
where the head of $e$ equals the tail of $e'$. Let $s$ be their common end and choose two disjoint paths $P_1,P_2$ of the torso of $s$ in $(T, (W_t:t\in V(T)))$ between the adhesion $\{u,v\}$
of $e$ and the adhesion $\{u',v'\}$ of $e'$, where $u,u'\in V(P_1)$. If $P_2$ has length at least one, then the model of $\F_{d_1}$ in $(G[A_{e'}], u',v')$ can be extended to show that 
$\F_{d_1+1}$ is a rooted minor of $(G[A_e],u,v)$, a contradiction; and similarly $P_1$ has length zero. But then the two adhesions are equal, contradicting that the tree-decomposition
is reduced. This proves \ref{fanthmagain}.~\bbox

\section{The proof of \ref{innerpthm}}
If $(G,u,v)$ is a 2-rooted graph, its {\em reverse} is $(G,v,u)$. 
We define $\IP_1 = \SP_1$, and for $d\ge 2$, even, $\IP_d$ denotes the right-legging of $\IP_{d-1}$. For $d\ge 2$, odd, $\IP_d$ is the parallel connection of $\IP_{d-1}$ and its reverse. (See Figure
\ref{fig:innerp}.)

\begin{figure}[h]
\begin{center}
\begin{tikzpicture}[scale = .9]
\tikzstyle{every node}=[inner sep=1.5pt, fill=black,circle,draw]
\node (1) at (1,0) {};
\node (2) at (2,0) {};
\node (3) at (3,0) {};
\node (5) at (4,0) {};
\node (6) at (5,0) {};
\node (9) at (6,0) {};

\node (10) at (7,0) {};
\node (14) at (8,0) {};
\node (15) at (9,0) {};
\node (19) at (11,0) {};
\tikzstyle{every node}=[inner sep=1.5pt, circle,draw]
\node (4) at (3.5,.6) {};
\node (7) at (5.2,.6) {};
\node (8) at (5.8,.6) {};
\node (11) at (7,.6) {};
\node (12) at (7.3,1) {};
\node (13) at (7.8,1) {};
\node (16) at (9,.6) {};
\node (17) at (9.3,1) {};
\node (18) at (9.8,1) {};
\node (20) at (11,.6) {};
\node (21) at (10.7,1) {};
\node (22) at (10.2,1) {};

\foreach \to/\from in {1/2,3/4,4/5,6/7,6/8,7/9,8/9,10/11,10/12,13/11,13/12,13/14,15/16,15/17,18/16,18/17,18/19,15/22,22/21,22/20,19/20,19/21}
\draw[-] (\to)--(\from);

\end{tikzpicture}
\end{center}
\caption{$\IP_1\LL \IP_5$.} \label{fig:innerp}
\end{figure}
Now we prove \ref{innerpthm}, which we restate:
\begin{thm}\label{innerpthmagain}
For every innerplanar graph $H$,
there exists $m$ such that every simple two-connected series-parallel graph with no $H$ minor admits a rooted tree-decomposition with adhesion two,
of height at most $m$, and such that each torso is either outerplanar or an apex-forest.
\end{thm}
\Proof
Every innerplanar graph is a minor of some $\IP_\ell$, so we may assume that $H=\IP_\ell$ for some $\ell$. Let $G$ be a two-connected simple series-parallel graph with 
no $\IP_\ell$ minor.
Choose a reduced 2-refined rooted tree-decomposition $(T, (W_t:t\in V(T)))$ as in \ref{goodtd}. For each edge $e\in E(T)$, let $\{u,v\}$ be the adhesion at $e$;
we say the score $\sc(e)$ of $e$ is $d_1+d_2$,
where $d_1$ is maximum such that $\IP_{d_1}$ is a rooted minor of $(G[A_e],u,v)$, and $d_2$ is maximum such that the reverse of $\IP_{d_2}$ is a rooted minor of $(G[A_e],u,v)$.
If $d>1$, then the reverse of $\IP_d$ contains $\IP_{d-1}$ as a rooted minor. Moreover, if $d$ is odd then $\IP_d$ is isomorphic to its own reverse. Consequently $d_1,d_2$ differ by at most one,
and if they are
different then the larger is even. 

Thus, all scores are between $2$ and $2\ell-2$; and as usual, for every path $P$ of $T$ with one end the root of $T$, the scores of the edges in $P$ monotonically decrease with distance from the
root. By \ref{genthm}, $G$ admits a rooted tree-decomposition with adhesion two and with height at most $4\ell$, such that each of its torsos is either 
 a triangle or the $S$-torso of some plateau $S$ with an internal vertex. It remains to show that for every plateau $S$ with an internal vertex, the $S$-torso $H(S)$ is
either an apex-forest, or outerplanar, and that is the content of the next result. This will prove \ref{innerpthmagain}.~\bbox

\begin{thm}\label{innerp} 
For every plateau $S$ with an internal vertex, the $S$-torso $H(S)$ is
either an apex-forest, or outerplanar. 
\end{thm}
\Proof
Let $e_0$ be the root-edge of $S$, let $d$ be the common score,  let $A_{e_0}\cap B_{e_0} = \{u_0,v_0\}$, and
choose $d_1,d_2$ as in the definition of score, with $d_1+d_2=d$. Let $r_0$ be the root of $S$. As in the proof of \ref{fanthm}, we have:
\\
\\
(1) {\em For every edge $e$ of $S$, let $P_1,P_2$ be disjoint paths between the adhesions of $e_0$ and of $e$, where $u_0\in V(P_1)$. Let the adhesion at $e$ be $\{u,v\}$, where $u\in V(P_1)$.
Then $d_1$ is maximum such that $\IP_{d_1}$ is a rooted minor of $(G[A_e],u,v)$, and $d_2$ is maximum such that the reverse of $\IP_{d_2}$ is a rooted minor of $(G[A_e],u,v)$.}

\bigskip

As we saw, either $d_1=d_2$, or they differ by one and the larger is even.
As in the proof of \ref{fanthm}, not both $d_1,d_2$ are odd, since $S$ has an internal vertex. Suppose they are both even, and hence equal. We proceed as in the proof of \ref{SPthm}, but we write
out the proof because there are slight differences.
Let $s\in V(S)\setminus \{r_0\}$, and let $r$ be its parent (thus $r\in V(S)$) and $f=rs$. Let $W_s=\{v_1,v_2,v_3\}$, where $A_f\cap B_f = \{v_1,v_2\}$. Since $(T,(W_t:t\in V(T)))$ is reduced, every edge
$e\ne f$
of $T$ incident with $s$ satisfies $A_e\cap B_e=\{v_1,v_3\}$ or  $A_e\cap B_e=\{v_2,v_3\}$.
\\
\\
(2) {\em There is at most one edge $e\in E(S)$ incident with $s$ with $A_e\cap B_e = \{v_1,v_3\}$, and the same for $\{v_2,v_3\}$.}
\\
\\
Since $d_1$ is even, $\IP_{d_1+1}$ is the parallel connection of $\IP_{d_1}$ and its reverse. 
Suppose that $e,e'\in E(T)$ are distinct edges of $S$ both incident with $s$ and both with adhesion
$\{v_1,v_3\}$. Since $e\in E(S)$, it follows that $\IP_{d_1}$ is a rooted minor of $(G[A_e],v_1,v_3)$. Since $d_2=d_1$, 
the reverse of $\IP_{d_1}$ is a rooted minor of $(G[A_{e'}],v_1,v_3)$. By combining the corresponding two models, it follows that the parallel
connection $\IP_{d_1+1}$ of $\IP_{d_1}$ and its reverse is a rooted minor of $(G[A_e\cup A_{e'}],v_1,v_3)$.
Let $\eta$ be the corresponding model of $\IP_{d_1+1}$.
Since $G$ is two-connected, there is a path of $G\setminus \{v_1\}$ between $v_2,v_3$, say $P$. By contracting this path, we convert $\eta$
to a model of $\IP_{d_1+1}$ in $G[A_f]$, and so $\IP_{d_1+1}$ is a rooted minor of $(G[A_f],v_1,v_2)$, a contradiction.
This proves (2).

\bigskip

Now every edge of the $S$-torso has both ends in $W_s$ for some $s\in V(S)\setminus \{r_0\}$. Let $J$ be the graph with vertex set $\bigcup_{s\in V(S)\setminus \{r_0\}}W_s$ and with edges every pair
$u,v$ of distinct vertices that both belong to some $W_s$ for some $s\in V(S)\setminus \{r_0\}$. From (2), $J$ is outerplanar, and hence so is the $S$-torso.

Thus we may assume that $d_1$ is even and $d_2=d_1-1$ is odd. Let the adhesion at $e_0$ be $\{u_0,v_0\}$, where $d_1$ is maximum such that $\IP_{d_1}$ is a rooted
minor of $(G[A_{e_0}],u_0,v_0)$.
\\
\\
(3) {\em For each $e\in E(S)$, the adhesion at $e$ contains $u_0$.}
\\
\\
Suppose that this is false for some $e\in E(S)$, with adhesion $\{u,v\}$ say. Choose two disjoint paths $P_1,P_2$ between the adhesions of $e$ and $e_0$, with
$u,u_0\in V(P_1)$ (exchanging the names of $u,v$ if necessary). Thus, $P_1$ has positive length. But $\IP_{d_2}$ (indeed, $\IP_{d_1}$, but all we need is $\IP_{d_2}$)
is a rooted minor of $(G[A_e], u,v)$. By adding $P_1$ to the corresponding model, we deduce that the left-legging of $\IP_{d_2}$  (that is, the reverse of
$\IP_{d_1}$) is a rooted minor of
$(G[A_{e_0}], u_0,v_0)$, contradicting that $d_2<d_1$. This proves (3).

\bigskip

From (3), we deduce that $H(S)$ is obtained from a graph of tree-width one by adding a vertex, and so is an apex-forest. This proves \ref{innerp}.~\bbox

\section{The proof of \ref{tw3}}

\begin{figure}[h]
\begin{center}
\begin{tikzpicture}[scale = .9]
\tikzstyle{every node}=[inner sep=1.5pt, fill=black,circle,draw]
\def\r{2}
\node (1) at ({\r*cos(22.5)},{\r*sin(22.5)}) {};
\node (2) at ({\r*cos(67.5)},{\r*sin(67.5)}) {};
\node (3) at ({\r*cos(112.5)},{\r*sin(112.5)}) {};
\node (4) at ({\r*cos(157.5)},{\r*sin(157.5)}) {};
\node (5) at ({\r*cos(202.5)},{\r*sin(202.5)}) {};
\node (6) at ({\r*cos(247.5)},{\r*sin(247.5)}) {};
\node (7) at ({\r*cos(292.5)},{\r*sin(292.5)}) {};
\node (8) at ({\r*cos(337.5)},{\r*sin(337.5)}) {};

\node (9) at (0,0) {};
\node (10) at (0,.7) {};
\node (11) at (-.6,1) {};
\node (12) at (.6,1) {};

\node (13) at (-.6,-.4) {};
\node (14) at (0,-.7) {};
\node (15) at (.6,-.4) {};

\foreach \to/\from in {1/2,2/3,3/4,4/5,5/6,6/7,7/8,8/1,9/10,10/11,10/12,9/13,9/14,9/15,11/5,11/3,11/4,12/2,12/1,15/1,15/7,14/7,14/6,13/6,13/5}
\draw[-] (\to)--(\from);

\end{tikzpicture}
\end{center}
\caption{A near-outerplanar graph.} \label{fig:nearouterp}
\end{figure}
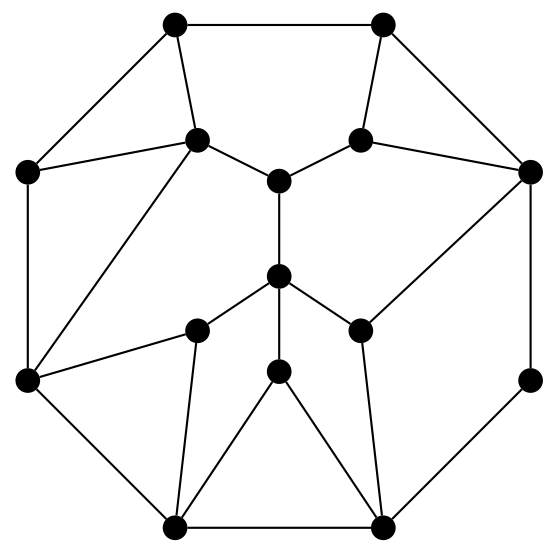

A graph is {\em near-outerplanar} if it can be drawn in the plane such that every region has a vertex in common with the infinite region. (See Figure \ref{fig:nearouterp}.)
If such a graph $G$ is simple and three-connected, the boundary of the infinite region is a cycle $C$, $G\setminus V(C)$ is a tree, and all its other edges are between this tree and $V(C)$. 

Near-outerplanar graphs have tree-width at most three, and form one of the building blocks of \ref{tw3}. Here is the full statement of \ref{tw3}:
\begin{thm}\label{tw3full}
For every series-parallel graph $H$ there exists $k$ such that every simple three-connected graph of tree-width three and with no $H$ minor admits a rooted tree-decomposition with adhesion three
and height at most $k$, such that each torso is either:
\begin{itemize}
\item a near-outerplanar graph; or
\item obtained from an outerplanar graph by adding a vertex; or
\item obtained from an innerplanar graph by adding a vertex; or
\item a graph with path-width at most three.
\end{itemize}
\end{thm}
\Proof
We may assume that $H$ is $\SP_h$. Let $G$ be a simple three-connected graph with tree-width three and  with no $H$ minor; then $G$ admits a rooted tree-decomposition of width three and adhesion three. 
As in the proof of \ref{goodtd}, we may choose a rooted tree-decomposition $(T,(W_t:t\in V(T)))$, of width three and adhesion three, such that 
\begin{itemize}
\item if $s,t\in V(T)$ are adjacent then $W_s\not\subseteq W_t$; and
\item for every path $P$ of $T$ with one end the root of $T$, no two edges of $P$ have the same adhesion. 
\end{itemize}
We need to define the score of edges of $T$. 
For each $e\in E(T)$ and all distinct $u,v\in A_e\cap B_e$, we define $\insc(u,v,e)$ to be the largest $d$ such that $\SP_d$ is a rooted minor of $(G',u,v)$, where
$G'=G[(A_e\setminus B_e)\cup \{u,v\}]$; and we define $\outsc(u,v,e)$ to be the largest $d$ such that $\SP_d$ is a rooted minor of $(G[A_e],u,v)$.
Thus, $\insc(u,v,e)\le \outsc(u,v,e)\le h-1$, since $G$ has no $\SP_h$ minor. 
We use this to define the score of edges of $T$, as follows.
Let $e\in E(T)$ with adhesion $\{u,v,w\}$; then 
$$\sc(e) = \insc(u,v,e) +\outsc(u,v,e) + \insc(u,w,e) +\outsc(u,w,e) + \insc(v,w,e) +\outsc(v,w,e)\le 6h.$$
We observe, first, that, as in the proof of \ref{SPthm}:
\\
\\
(1) {\em Let $e,e'$ be distinct edges of $T$ where $e$ belongs to the path between $e'$ and the root of $T$; and choose three disjoint paths $P_1,P_2,P_3$ between 
the corresponding two adhesions.
For $i = 1,2,3$, let $v_i,v_i'$ be the vertices of $P_i$ that belong to the adhesion at $e$ and the adhesion at $e'$, respectively. Then $\insc(v_i,v_j,e)\ge \insc(v_i',v_j',e')$
and $\outsc(v_i,v_j,e)\ge \outsc(v_i',v_j',e')$ for all distinct $i,j\in \{1,2,3\}$. Hence $\sc(e')\le \sc(e)$, and if $\sc(e')= \sc(e)$, then 
$\insc(v_i,v_j,e)= \insc(v_i',v_j',e')$ and $\outsc(v_i,v_j,e)= \outsc(v_i',v_j',e')$ for all distinct $i,j\in \{1,2,3\}$.}

\bigskip

We define plateaus as before; and for each plateau $S$ with an internal vertex, we define the $S$-torso $H(S)$ as follows.
Let $G_1$ be the union of $G[W_s]$ over
all $s\in V(S^-)$. The 
$S$-torso is obtained from $G_1$ by adding new edges joining every two of the three members of $A_e\cap B_e$ (if they are currently nonadjacent), for each edge $e\in E(T)$ that has exactly one end in $V(S^-)$.
We deduce as before that 
$G$ admits a rooted tree-decomposition with adhesion three and of height at most $12k$, such that each of its torsos
is either $H(S)$ for some plateau $S$ with an internal vertex, or is a complete graph $K_4$. We need to figure out the possibilities for the graphs $H(S)$, and they are handled in
the following results. Let $e_0$ be the root-edge of $S$, and let the adhesion at $e_0$ be $\{u_1,u_2,u_3\}$. For $h = 1,2,3$, define 
\begin{align*}
c_h &=\outsc(u_i,u_j,e_0)\\
d_h &= \insc(u_i,u_j, e_0)
\end{align*}
where $\{i,j\} = \{1,2,3\}\setminus \{h\}$.

We need to show that $H(S)$ is one of the four bullets of the theorem. 
\\
\\
(2) {\em If $c_1\ne d_1$ then $H(S)$ is as in one of the four bullets.}
\\
\\
Suppose that $c_1\ne d_1$. Hence $\SP_{c_1}$ is contained as a rooted minor of $(G[A_{e_0}],u_2,u_3)$, but not of $(G[A_{e_0}\setminus \{u_1\}],u_2,u_3)$.
Suppose that for some $e\in E(S)$, the adhesion at $e$ does not contain $u_1$. Let the adhesion at $e$ be $\{v_1,v_2,v_3\}$ where there are disjoint paths $P_1,P_2,P_3$
such that $P_i$ links $u_i$ with 
$v_i$ for $i = 1,2,3$. Since $\outsc(v_2,v_3,e) = c_1$, $\SP_{c_1}$ is contained as a rooted minor of $(G[A_{e}],v_2,v_3)$. But since $u_1\notin A_e$, and $u_1,v_1$ do not 
belong to $P_2,P_3$, it follows that 
$\SP_{c_1}$ is contained as a rooted minor of $(G[A_{e_0}\setminus \{u_1\}],u_2,u_3)$, contradicting that $d_1<c_1$. 
So $u_1$ belongs to all the adhesions at $e$ (in  $(T,(W_t:t\in V(T)))$) for every edge $e$ of $S$. But then $H(S)\setminus \{u_1\}$ is either outerplanar or innerplanar, 
as in the proof of \ref{SPthm}, depending on the parity of $d_1$, and so $H(S)$ is indeed as in one of the four bullets. This proves (2).

\bigskip

So we may assume that $c_i = d_i$ for $i = 1,2,3$.  Henceforth we write $\sc$ for both $\insc$ and $\outsc$. Let $r_0$ be the root of $S$. 
\\
\\
(3) {\em Let $s\in V(S)$. If $s_1,s_2\in V(S)$ are two of its children, then the adhesions of $ss_1, ss_2$ are different.}
\\
\\
Since $s$ has more than one child in $S$, it follows that $s\ne r_0$; let $r$ be the parent of $s$. Let 
$W_s = \{v_1,v_2,v_3,w\}$ where $\{v_1,v_2,v_3\}$ is the adhesion at $rs$. 
We may assume that  $\sc(v_1,v_2,rs)= d_3$, and $\sc(v_1,v_3,rs) = d_2$, and $\sc(v_2,v_3,rs) = d_1$.
The adhesion at $ss_1$ is a subset of $W_s$ of size three, and not $\{v_1,v_2,v_3\}$
from the choice of $(T,(W_t:t\in V(T)))$. Suppose that the adhesions of $ss_1$ and $ss_2$ are equal; then we may assume they both equal $\{v_1,v_2,w\}$.
Since $G$ is three-connected, there is a path $P_3$ between $w, v_3$ not containing $v_1$ or $v_2$. From $P_3$ and the two one-vertex paths with vertices $v_1,v_2$ respectively,
we deduce from (1) that for $i = 1,2$, $\sc(v_1,v_2,ss_i) = d_3$, and $\sc(v_1,w,ss_i) = d_2$, and $\sc(v_2,w,ss_i) = d_1$. 

There is a model of $\SP_{d_3}$ in $(G[A_{ss_1}\setminus \{w\}], v_1,v_2)$, and there is a path $Q$ of $G[A_{ss_2}\setminus \{v_2\}]$ between $v_1,w$. By contracting $Q$
and $P_3$,
we obtain a model of $\SP_{d_3}$ in $(G[A_{rs}],v_2,v_3)$, and we deduce that $\sc(v_2,v_3,rs)\ge d_3$.  Since $d_1 = c_1 = \sc(v_2,v_3,rs)$ it follows that $d_1\ge d_3$. 
Similarly, $d_2\ge d_3$. Now there are two cases. Suppose that $d_3$ is even;  then $\SP_{d_3+1}$ is the parallel doubling of $\SP_{d_3}$. By combining the models of 
$\SP_{d_3}$ in $(G[A_{ss_i}\setminus \{w\}],v_1,v_2)$ for $i = 1,2$, we obtain a model of the parallel doubling of $\SP_{d_3}$ (that is, $\SP_{d_3+1}$) in
$(G[A_{rs}],v_1,v_2)$, contradicting that $\outsc(v_1,v_2,rs) = d_3$. Now we suppose that $d_3$ is odd, and so $\SP_{d_3+1}$ is the series doubling of $\SP_{d_3}$.
Since $d_1,d_2\ge d_3$, there are models of $\SP_{d_3}$ in $(G[A_{ss_1}\setminus \{v_1\}],v_2,w)$ and in 
$(G[A_{ss_2}\setminus \{v_2\}],v_1,w)$. By combining them, we obtain a model of the series doubling of $\SP_{d_3}$ (that is, $\SP_{d_3+1}$) in $(G[A_{rs}], v_1,v_2)$,
contradicting that $\sc(v_1,v_2,rs) = d_3$. This proves (3).

\bigskip

\begin{figure}[h]
\begin{center}
\begin{tikzpicture}[scale = 1.2]
\tikzstyle{every node}=[inner sep=1.5pt, fill=black,circle,draw]
\node (s) at (-4,.5) {};
\node (s2) at (-5.2,1.5) {};
\node (s1) at (-2.8,1.5) {};
\node (r) at (-4,-1) {};
\foreach \to/\from in {s/s1,s/s2,s/r}
\draw[-] (\to)--(\from);

\node (v1) at (-2,0) {};
\node (v2) at (2,0) {};
\node (v3) at (0,0) {};
\node (w) at (0,2) {};
\draw[fill = yellow] (v1) to [bend left = 20] (v3) to [bend left = 20] (w) to [bend right = 20] (v1);
\draw[fill = green] (v3) to [bend left = 20] (v2) to [bend right = 20] (w) to [bend left = 20] (v3);
\node (v3) at (0,0) {};
\draw[-] (v1) to (-2.1,-.4);
\draw[-] (v1) to (-1.9,-.4);
\draw[-] (v2) to (2.1,-.4);
\draw[-] (v2) to (1.9,-.4);
\draw[-] (v3) to (0.1,-.4);
\draw[-] (v3) to (-.1,-.4);

\tikzstyle{every node}=[]
\draw (s) node[left] {$s$};
\draw (s1) node[right] {$s_1$};
\draw (s2) node[left] {$s_2$};
\draw (r) node[left] {$r$};
\draw (v1) node[left] {$v_1$};
\draw (v2) node[right] {$v_2$};
\draw (w) node[above] {$w$};
\draw (v3) node[below right] {$v_3$};
\node at (-1,.7) {$A_{ss_2}$};
\node at (1,.7) {$A_{ss_1}$};

\end{tikzpicture}
\end{center}
\caption{Two children with different adhesions.} \label{fig:boat}
\end{figure}
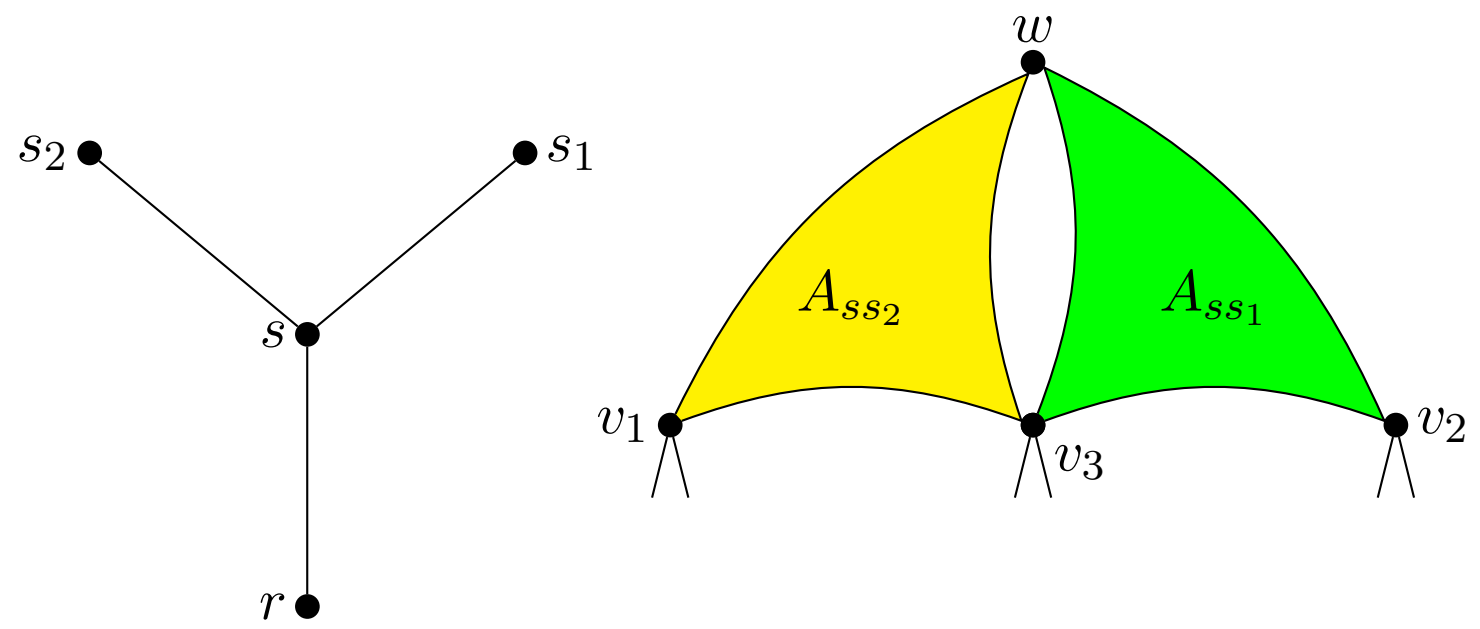

If $S$ is a path, then $H(S)$ has path-width at most three and satisfies one of the bullets of the theorem, so henceforth we assume that $S$ is not a path.
\\
\\
(4) {\em $d_1=d_2=d_3=d$ say, and $d$ is even.}
\\
\\
Since $S$ is not a path, there exists 
$s\in V(S)$, with children $s_1,s_2\in V(S)$. Let $r$ be the parent of $s$, and let $W_s = \{v_1,v_2,v_3,w\}$ where $\{v_1,v_2,v_3\}$ is the adhesion at $rs$.
By (3), the adhesions of $ss_1, ss_2$ are different, so we may assume that the adhesion at $ss_i$ is
$W_s\setminus \{v_i\}$ for $i = 1,2$. (See Figure \ref{fig:boat}.)

We may assume that  $\sc(v_1,v_2,rs) = d_3$, and $\sc(v_1,v_3,rs) = d_2$, and $\sc(v_2,v_3,rs) = d_1$.
As in (3), $\sc(v_2,v_3,ss_1) = d_1$, and $\sc(v_2,w,ss_1) = d_3$, and
$\sc(v_3,w,ss_1) = d_2$, and similarly $\sc(v_1,v_3,ss_2) = d_2$, and $\sc(v_1,w,ss_2) = d_3$, and $\sc(v_3,w,ss_2) = d_1$.
By contracting  a path in $G[A_{ss_2}\setminus \{w\}]$ between $v_1,v_3$, we convert a model of $\SP_{d_1}$ in $(G[A_{ss_1}\setminus \{w\}], v_2,v_3)$ 
to a model of $\SP_{d_1}$ in $(G[A_{rs}], v_1,v_2)$, and so $\sc(v_1,v_2,rs)\ge d_1$; and consequently $d_3\ge d_1$, and similarly $d_3\ge d_2$. But by 
contracting a path in $G[A_{ss_2}\setminus \{v_1\}]$ between $w,v_3$, we convert 
 a model of $\SP_{d_3}$ in $(G[A_{ss_1}\setminus \{v_3\}], v_2,w)$ to 
a model of $\SP_{d_3}$ in $(G[A_{rs}], v_2,v_3)$, and so $\sc(v_2,v_3,rs)\ge d_3$; and therefore $d_1\ge d_3$, and similarly $d_2\ge d_3$. So $d_1=d_2=d_3=d$ say.

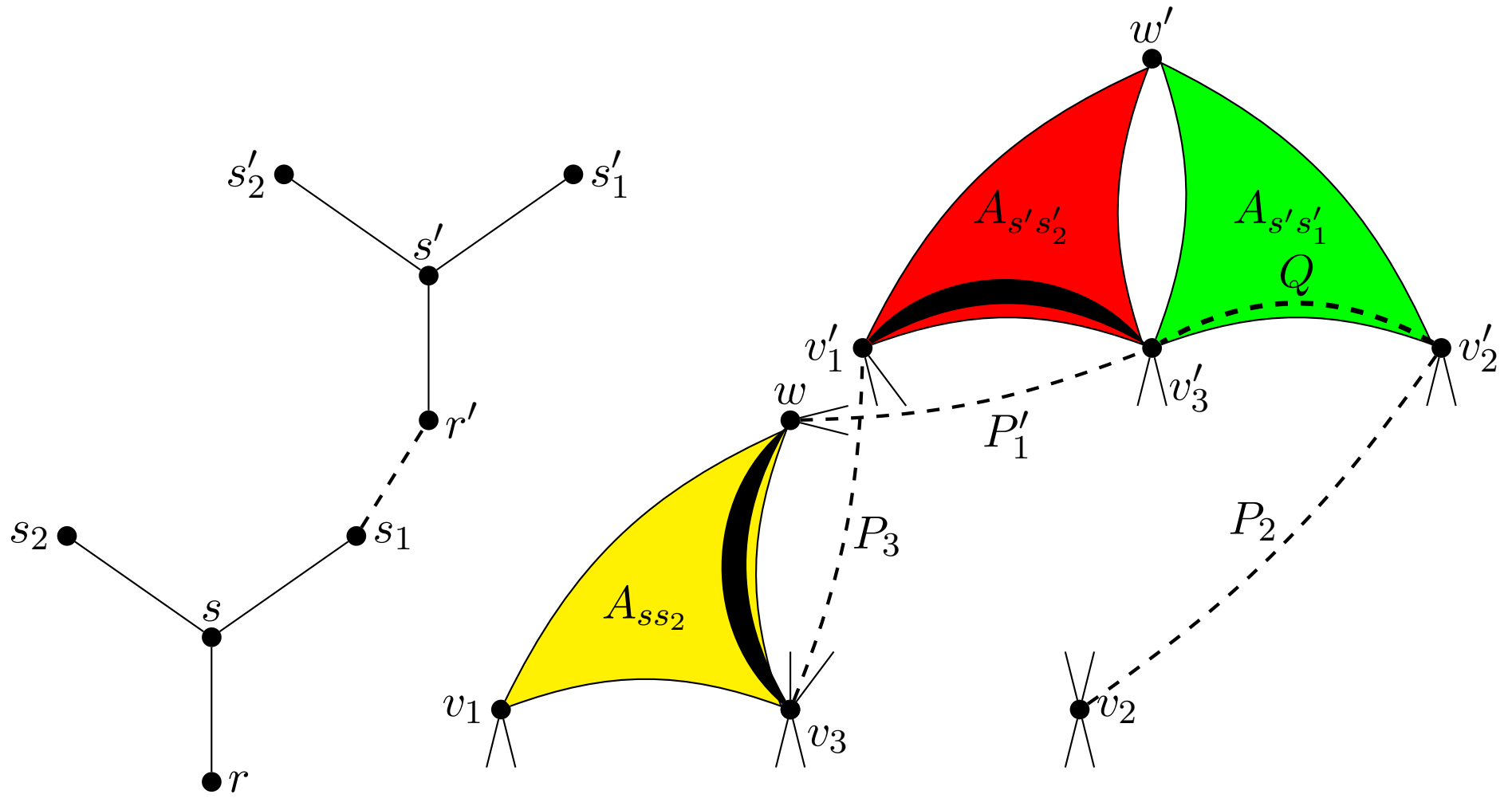
\begin{figure}[h]
\begin{center}
\begin{tikzpicture}[scale = 1.2]
\tikzstyle{every node}=[inner sep=1.5pt, fill=black,circle,draw]
\node (s) at (-4,.5) {};
\node (s2) at (-5,1.2) {};
\node (s1) at (-3,1.2) {};
\node (r) at (-4,-.5) {};
\node (r') at (-2.5,2) {};
\node (s') at (-2.5, 3) {};
\node (s1') at  (-1.5, 3.7) {};
\node (s2') at  (-3.5, 3.7) {};
\foreach \to/\from in {s/s1,s/s2,s/r, s'/s1',s'/s2',s'/r'}
\draw[-] (\to)--(\from);
\draw[dashed, thick] (s1) to (r');

\node (v1) at (-2,0) {};
\node (v2) at (2,0) {};
\node (v3) at (0,0) {};
\node (w) at (0,2) {};
\node (v1') at (.5,2.5) {};
\node (v2') at (4.5,2.5) {};
\node (v3') at (2.5,2.5) {};
\node (w') at (2.5,4.5) {};

\draw[fill = yellow] (v1) to [bend left = 20] (v3) to [bend left = 20] (w) to [bend right = 20] (v1);
\draw[fill = red] (v1') to [bend left = 20] (v3') to [bend left = 20] (w') to [bend right = 20] (v1');
\draw[fill = green] (v3') to [bend left = 20] (v2') to [bend right = 20] (w') to [bend left = 20] (v3');

\draw[fill=black] (v3) to [bend left = 30] (w) to [bend right = 50] (v3);
\draw[fill=black] (v1') to [bend left = 30] (v3') to [bend right = 50] (v1');
\node (v3) at (0,0) {};
\draw[-] (v1) to (-2.1,-.4);
\draw[-] (v1) to (-1.9,-.4);
\draw[-] (v2) to (2.1,-.4);
\draw[-] (v2) to (1.9,-.4);
\draw[-] (v3) to (0.1,-.4);
\draw[-] (v3) to (-.1,-.4);

\draw[dashed, thick] (v3) to [bend right = 10] (v1');
\draw[dashed, thick] (w) to [bend right = 10] (v3');
\draw[dashed, thick] (v2) to [bend right = 10] (v2');
\draw[dashed, very thick] (v2') to [bend right = 30] (v3');

\tikzstyle{every node}=[]
\draw (s) node[above] {$s$};
\draw (s1) node[right] {$s_1$};
\draw (s2) node[left] {$s_2$};
\draw (r) node[right] {$r$};
\draw (s') node[above] {$s'$};
\draw (s1') node[right] {$s_1'$};
\draw (s2') node[left] {$s_2'$};
\draw (r') node[right] {$r'$};

\draw (v1) node[left] {$v_1$};
\draw (v2) node[right] {$v_2$};
\draw (w) node[above] {$w$};
\draw (v3) node[below right] {$v_3$};
\node at (-1,.7) {$A_{ss_2}$};
\node at (1.6,3.4) {$A_{s's_2'}$};
\node at (3.4,3.4) {$A_{s's_1'}$};

\draw (v1') node[left] {$v_1'$};
\draw (v2') node[right] {$v_2'$};
\draw (w') node[above] {$w'$};
\draw (v3') node[below right] {$v_3'$};
\node at (1.5,1.9) {$P_1'$};
\node at (3.2,1.3) {$P_2$}; 
\node at (.6,1.2) {$P_3$};
\node at (3.5,3) {$Q$};

\draw[-] (v3') to (2.6,2.1);
\draw[-] (v3') to (2.4,2.1);
\draw[-] (v2') to (4.6,2.1);
\draw[-] (v2') to (4.4,2.1);
\draw[-] (v1') to (.6,2.1);
\draw[-] (v1') to (.8,2.1);
\draw[-] (v2) to (2.1,.4);
\draw[-] (v2) to (1.9,.4);
\draw[-] (v3) to (.3,.4);
\draw[-] (v3) to (0,.4);
\draw[-] (w) to (.4,1.9);
\draw[-] (w) to (.4,2.1);

\end{tikzpicture}
\end{center}
\caption{For the proof of (6). The black bananas represent models of $\SP_d$. The paths $P_1'$ and $P_2$ might link $w$ to $v_2'$ and $v_2$ to $v_3'$ instead of what is drawn.} \label{fig:straightpath}
\end{figure}

By combining a model of $\SP_d$ in $(G[A_{ss_2}\setminus \{v_3\}], v_1,w)$ with a model of $\SP_d$ in $(G[A_{ss_1}\setminus \{v_3\}], v_2,w)$, we obtain a model 
of the series doubling of $\SP_d$ in $(G[A_{rs}],v_1,v_2)$. Since there is no model of $\SP_{d+1}$ in $(G[A_{rs}],v_1,v_2)$, we deduce that $d$ is even. This proves (4).
\\
\\
(5) {\em Every vertex in $S$ has at most two children in $S$.}
\\
\\
Suppose that $s\in V(S)$ has three children $s_1,s_2,s_3$ in $S$; and let $r$ be the parent of $s$ in $S$, and $W_s = \{v_1,v_2,v_3,w\}$. By (3), we may assume that 
the adhesion at $ss_i$ is $W_s\setminus \{v_i\}$ for $i = 1,2,3$.
By combining the models of $\SP_d$ in $(G[A_{ss_1}\setminus \{v_2\}],v_3,w)$ and in $(G[A_{ss_2}\setminus \{v_1\}],v_3,w)$ we obtain a model of the parallel doubling of $\SP_d$, 
that is, of $\SP_{d+1}$, in 
$$(G[(A_{ss_1}\setminus \{v_2\})\cup (A_{ss_2}\setminus \{v_1\})],v_3,w)=(G[(A_{ss_1}\cup A_{ss_2})\setminus \{v_1,v_2\}],v_3,w).$$
By contracting a path of $G[A_{ss_3}\setminus \{v_2\}]$ between $w,v_1$, we convert this model to a model of $\SP_{d+1}$ in $(G[A_{rs}],v_1,v_3)$, a contradiction.
This proves (5). 
\\
\\
(6) {\em 
Let $s,s'\in V(S)$, with parents $r,r'$ respectively, where $s$ is in the path of $S$ between $e_0$ and $s'$. Let $s$ have two children $s_1,s_2$ in $S$, and let $s'$ have two children $s_1',s_2'$ in  $S$. 
Let $W_s=\{v_1,v_2,v_3,w\}$, where the adhesion at $rs$ is $\{v_1,v_2,v_3\}$, and the adhesion of $ss_i$ is $W_s\setminus \{v_i\}$ for $i= 1,2$; 
and let $W_{s'} = \{v_1',v_2',v_3',w'\}$, where the adhesion at $r's'$ is $\{v_1',v_2',v_3'\}$, and the adhesion 
of $s's_i'$ is $W_{s'}\setminus \{v_i'\}$ for $i= 1,2$. For every set of three vertex-disjoint paths between $\{v_1,v_2,v_3\}$
and $\{v_1',v_2',v_3'\}$, one of the three paths joins $v_3$ and $v_3'$.}
\\
\\
We may assume that $s\ne s'$. 
Since $s'$ is a descendant of $s$ and not equal to $s$, it is a descendant of one of $s_1,s_2$ and we assume $s'$
is a descendant of $s_1$ without loss of generality. Suppose that $P_1,P_2,P_3$ are vertex-disjoint paths between $\{v_1,v_2,v_3\}$
and $\{v_1',v_2',v_3'\}$, and none of them joins $v_3,v_3'$. Let $v_i\in V(P_i)$ for $i = 1,2,3$.
So $w\in V(P_1)$, since each of $P_1,P_2,P_3$ has a vertex in $\{w,v_2,v_3\}$. Let $P_1'$ be the subpath of $P_1$ between $w$ and the end of $P_1$ in $\{v_1',v_2',v_3'\}$. 
(So any or all of $P_1', P_2,P_3$ might have length zero.)
Since $P_3$ does not contain $v_3'$, we may assume
(by exchanging $v_1',v_2'$ if necessary) that $v_1'\in V(P_3)$; and then there are two cases: $v_3'\in V(P_1)$ and $v_2'\in V(P_2)$; and $v_3'\in V(P_2)$ and $v_2'\in V(P_1)$.
(See Figure \ref{fig:straightpath}.)
In either case, let $Q$ be a path in $G[A_{s's_1'}\setminus \{w'\}]$ between $v_2',v_3'$; so $P_1'\cup Q\cup P_2$ is a path between $w, v_2$ passing through $v_3'$ and disjoint from $P_3$.
By contracting $P_1'\cup Q\cup P_2$ and $P_3$, 
we convert the model of $\SP_d$ in $(G[A_{ss_2}\setminus \{v_1\}], v_3,w)$ and the model of $\SP_d$ in $(G[A_{s's_2'}\setminus \{w'\}],v_1',v_3')$
to a model of $\SP_{d+1}$ in  $(G[A_{rs}], v_2,v_3)$, which is impossible. This proves (6).

\bigskip
Let us define (inductively on the length of the path between $e_0,e$) a numbering, called the ``standard numbering'', of the three members of the adhesion at $e$, for each edge $e\in E(S)$, as follows. 
The adhesion at $e_0$ has standard numbering $(u_1,u_2,u_3)$. For each edge $e\in E(S)$ with $e\ne e_0$, let $s$ be the end of $e$ in the interior of the path between $e,e_0$, 
let $r$ be the parent of $s$, and let $e'=rs$. We assume inductively that the standard numbering of the adhesion at $e'$ is defined, say $(v_1,v_2,v_3)$. Since
the adhesions at $e$ and at $e'$ are both subsets of $W_s$, and $|W_s|=4$, the adhesion of $e$ contains two of $v_1,v_2,v_3$, say the two different from $v_h$; and its third member is the 
vertex of $W_s$ different from $v_1,v_2,v_3$, say $w$. Then the standard numbering of the adhesion at $e$ is defined to be the sequence obtained from $(v_1,v_2,v_3)$
by replacing $v_h$ by $w$.

It follows immediately that if $e\in E(S)$ and $(v_1,v_2,v_3)$ is the standard numbering of the adhesion at $e$, and 
$P_1,P_2,P_3$ are pairwise vertex-disjoint paths of $G$ (or of $H(S)$) between the adhesions of $e_0$ and $e$, and $u_i\in V(P_i)$ for $i = 1,2,3$,
then $v_i\in V(P_i)$ for $i = 1,2,3$. 

If $s\in V(S)$ with parent $r$ and two children in $S$, we say that ``$s$ divides according to the standard numbering'' if  we can number the children $s_1,s_2$
such that the adhesion of $ss_i$ is $W_s\setminus \{v_i\}$ for $i = 1,2$,
where $W_s = \{v_1,v_2,v_3,w\}$, and $(v_1,v_2,v_3)$ is the standard numbering of the adhesion at $rs$. 
\\
\\
(7) {\em We may assume that every vertex of $S$ with two children
divides according to the standard numbering.}
\\
\\
Choose $s\in V(S)$ minimal (that is, as close to $r_0$ as possible) with 
two children. We may then choose the standard numbering such that $s$ divides according to the standard numbering. From (6), it follows that every vertex of $S$ with two children
also divides according to the standard numbering. This proves (7).
\\
\\
(8) {\em For every edge $e=rs$ of $S$, if
$(v_1,v_2,v_3)$ is
the standard numbering of the adhesion at $e$, then $v_1v_2\notin E(G)$, and if $e\ne e_0$, there is no child $t$ of $r$ in $T$ with $t\ne s$ such that $\{v_1,v_2\}\subseteq W_t$.}
\\
\\
Suppose this is false; then there is a path $Q$ of $G$ between $v_1,v_2$, either of length one, or contained in $A_{rt}$ for some child $t$ of $r$ different from $s$
with $v_1,v_2\in W_t$, and in this case $e\ne e_0$. If $t$ exists then $t\notin V(S)$, by (7).

Since $S$ 
is not a path, there is a vertex in $S$ with two children in $S$; choose such a vertex, $s'$ say, with minimum distance from $r_0$. 
Let $s_1',s_2'$ be its two children in $S$, $r'$ be its parent, and let $(v_1',v_2',v_3')$ be the standard numbering of the adhesion at $r's'$. Let $w'$ be the fourth vertex of $W_{s'}$.

\begin{figure}[h]
\begin{center}
\begin{tikzpicture}[scale = 1.2]
\tikzstyle{every node}=[inner sep=1.5pt, fill=black,circle,draw]
\node (s) at (-4,.5) {};
\node (r) at (-4,-.5) {};
\node (r') at (-4,2) {};
\node (s') at (-4, 3) {};
\node (s1') at  (-3, 3.7) {};
\node (s2') at  (-5, 3.7) {};
\node (q) at (-4,-1.5) {};
\node (t) at (-4.5,.5) {};

\foreach \to/\from in {q/r,r/t,s/r, s'/s1',s'/s2',s'/r'}
\draw[-] (\to)--(\from);
\draw[dashed, thick] (s) to (r');

\node (v1) at (-2,0) {};
\node (v2) at (2,0) {};
\node (v3) at (0,0) {};
\node (v1') at (-2,2.5) {};
\node (v2') at (2,2.5) {};
\node (v3') at (0,2.5) {};
\node (w') at (-0,4.5) {};

\draw[fill = yellow] (v1') to [bend left = 20] (v3') to [bend left = 20] (w') to [bend right = 20] (v1');
\draw[fill = green] (v3') to [bend left = 20] (v2') to [bend right = 20] (w') to [bend left = 20] (v3');

\draw[fill=black] (v1') to [bend left = 30] (v3') to [bend right = 50] (v1');
\draw[fill=black] (v3') to [bend left = 30] (v2') to [bend right = 50] (v3');

\draw[dashed, thick] (v3) to (v3');
\draw[dashed, thick] (v1) to (v1');
\draw[dashed, thick] (v2) to (v2');

\draw[dashed, thick] (v1) to [bend left = 20] (v2);

\tikzstyle{every node}=[]
\draw (s) node[right] {$s$};
\draw (r) node[right] {$r$};
\draw (s') node[above] {$s'$};
\draw (s1') node[right] {$s_1'$};
\draw (s2') node[left] {$s_2'$};
\draw (r') node[right] {$r'$};
\draw (t) node[left] {$t$};
\node at (-4.2,-1) {$g$};

\draw (v1) node[left] {$v_1$};
\draw (v2) node[right] {$v_2$};
\draw (v3) node[right] {$v_3$};

\node at (-.9,3.4) {$A_{s's_2'}$};
\node at (.9,3.4) {$A_{s's_1'}$};

\draw (v1') node[left] {$v_1'$};
\draw (v2') node[right] {$v_2'$};
\draw (w') node[above] {$w'$};
\draw (v3') node[below right] {$v_3'$};
\node at (-2.2,1.2) {$P_1$};
\node at (2.2,1.2) {$P_2$};
\node at (.2,1.2) {$P_3$};
\node at (-.6, .6) {$Q$};

\end{tikzpicture}
\end{center}
\caption{For the proof of (8).}
\label{fig:step9a}
\end{figure}
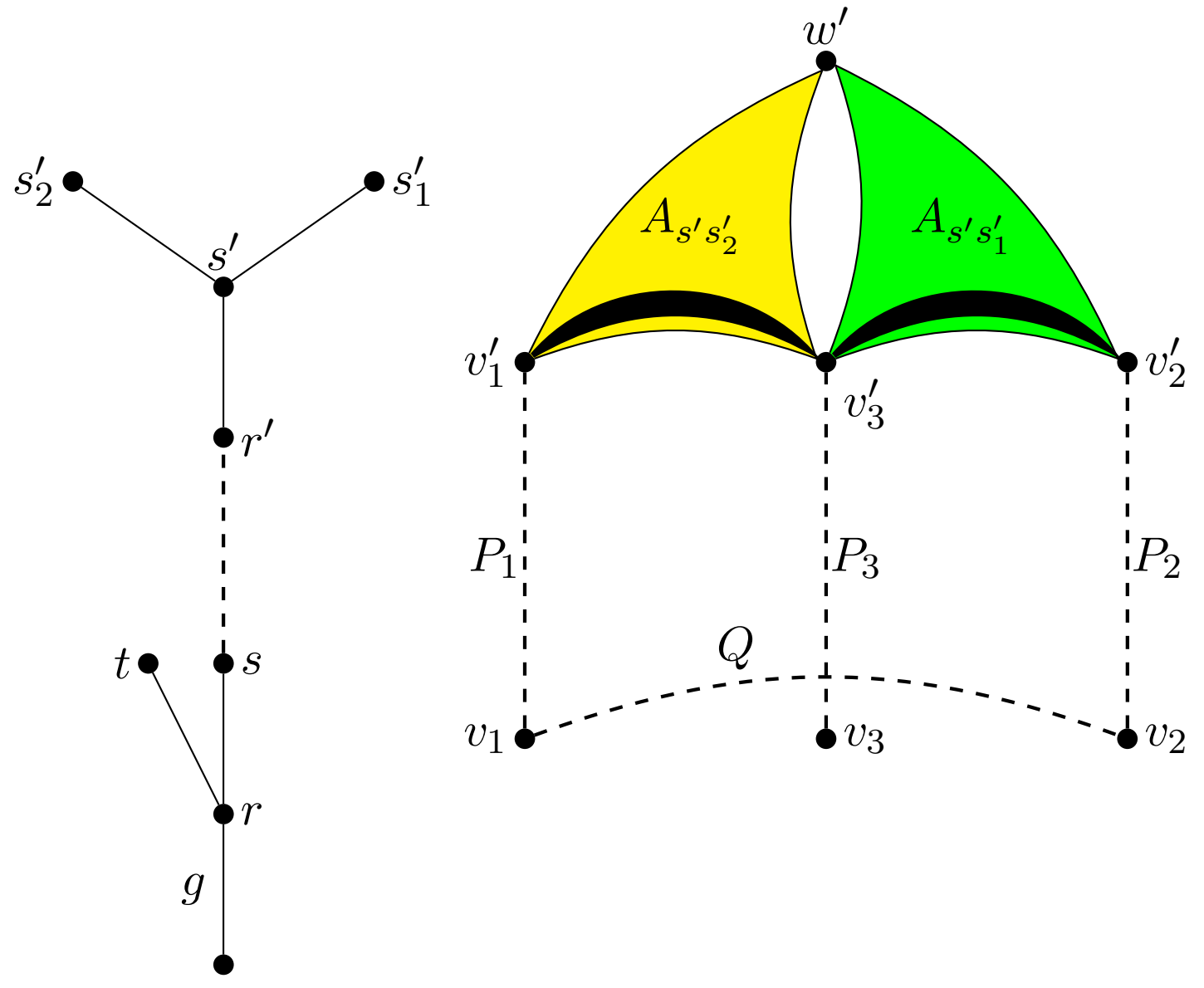

Now $s'$ is either an ancestor or a
descendant of $s$.
Suppose first that 
$s'$ is a descendant of $s$, possibly equal to $s$. Choose three disjoint paths $P_1,P_2,P_3$ of $G$ between the adhesions of $rs$ and $r's'$,
where $P_i$ joins $v_i, v_i'$ for $i = 1,2,3$.  (See Figure \ref{fig:step9a}.) 
If $Q$ has length one, then by contracting  $P_1,P_2, P_3$ and the edge $v_1v_2$, we convert the model of $\SP_d$
in $(G[A_{s's_2'}\setminus \{w'\}], v_1',v_3')$ and the model  of $\SP_d$
in $(G[A_{s's_1'}\setminus \{w'\}], v_2',v_3')$ to a model of $\SP_{d+1}$ in $(G[A_{rs}],v_2,v_3)$, a contradiction.
So we may assume that $Q$ has length more than one, and so $t$ exists, and $e\ne e_0$. Let $g$ be the edge of $S$ between $r$ and its parent, and let $W_r = \{v_1,v_2,v_3,x\}$ say, where
the adhesion at $g$ is $W_r\setminus \{v_i\}$ for some $i\in \{1,2,3\}$.  The adhesion at $rt$ contains all of $v_i, v_1,v_2$, but $v_i$ might be equal to one of $v_1,v_2$; so 
either
\begin{itemize}
\item $i=3$, the adhesion at $rt$ is $\{v_1,v_2,v_3\}$, and the adhesion at $g$ is $\{v_1,v_2,x\}$; or
\item $i\in \{1,2\}$, say $i = 1$, and the adhesion at $rt$ is $\{v_1,v_2,v_3\}$, and the adhesion at $g$ is $\{v_2,v_3,x\}$; or
\item $i\in \{1,2\}$, say $i = 1$, and the adhesion at $rt$ is $\{v_1,v_2,x\}$, and the adhesion at $g$ is $\{v_2,v_3,x\}$.
\end{itemize}
Now $A_{rs}\setminus W_r$ is vertex-disjoint from $A_{rt}\setminus W_r$; so 
the interior of $Q$ is disjoint from $A_{s's_1'}, A_{s's_2'}$ and the interior of $P_3$. 
If $v_2,v_3$ both belong to the adhesion at $g$, then
by contracting  $P_1,P_2, P_3, Q$, we convert the model of $\SP_d$
in $(G[A_{s's_2'}\setminus \{w'\}], v_1',v_3')$ and the model  of $\SP_d$
in $(G[A_{s's_1'}\setminus \{w'\}], v_2',v_3')$ to a model of $\SP_{d+1}$ in $(G[A_{g}],v_2,v_3)$, a contradiction.
So we assume that the adhesion at $rt$ is $\{v_1,v_2,v_3\}$, and the adhesion at $g$ is $\{v_1,v_2,x\}$. Since $G$ is three-connected, there is a path $R$ of $G\setminus \{v_1,v_2\}$ 
between $x,v_3$, and then by contracting  $P_1,P_2, P_3, Q, R$, we convert the model of $\SP_d$
in $(G[A_{s's_2'}\setminus \{w'\}], v_1',v_3')$ and the model  of $\SP_d$
in $(G[A_{s's_1'}\setminus \{w'\}], v_2',v_3')$ to a model of $\SP_{d+1}$ in $(G[A_{g}],v_2,x)$, a contradiction.

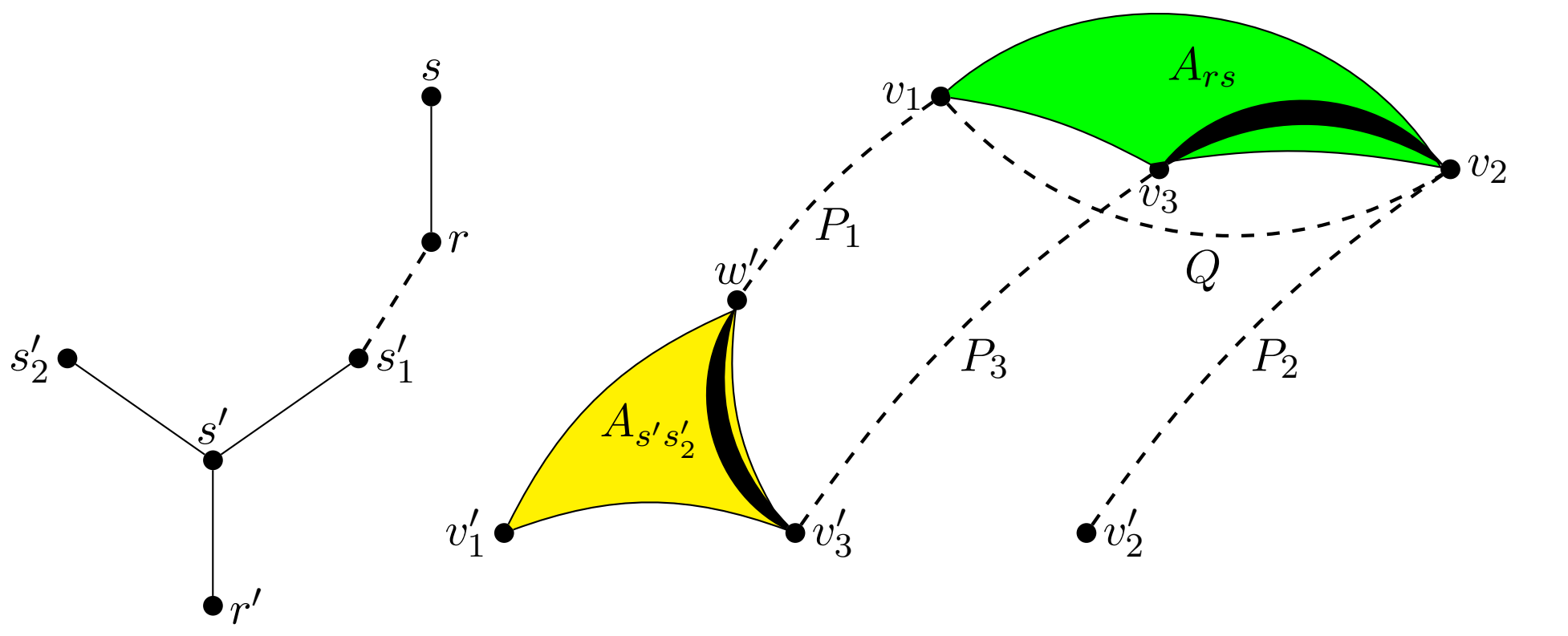
\begin{figure}[h]
\begin{center}
\begin{tikzpicture}[scale = 1.2]
\tikzstyle{every node}=[inner sep=1.5pt, fill=black,circle,draw]
\node (s') at (-4,.5) {};
\node (s2') at (-5,1.2) {};
\node (s1') at (-3,1.2) {};
\node (r') at (-4,-.5) {};
\node (r) at (-2.5,2) {};
\node (s) at (-2.5, 3) {};
\foreach \to/\from in {s'/s1',s'/s2',s'/r', s/r}
\draw[-] (\to)--(\from);
\draw[dashed, thick] (s1') to (r);

\node (v1') at (-2,0) {};
\node (v2') at (2,0) {};
\node (v3') at (0,0) {};
\node (w') at (-.4,1.6) {};
\node (v1) at (1,3) {};
\node (v2) at (4.5,2.5) {};
\node (v3) at (2.5,2.5) {};

\draw[fill = yellow] (v1') to [bend left = 20] (v3') to [bend left = 20] (w') to [bend right = 20] (v1');
\draw[fill = green] (v1) to [bend left = 10] (v3) to [bend left = 10] (v2) to [bend right = 50] (v1);

\draw[fill=black] (v3') to [bend left = 30] (w') to [bend right = 50] (v3');
\draw[fill=black] (v3) to [bend left = 30] (v2) to [bend right = 50] (v3);

\draw[dashed, thick] (v3) to [bend right = 10] (v3');
\draw[dashed, thick] (v1) to [bend right = 10] (w');
\draw[dashed, thick] (v2) to [bend right = 10] (v2');

\draw[dashed, thick] (v1) to [bend right = 40] (v2);

\tikzstyle{every node}=[]
\draw (s) node[above] {$s$};
\draw (s1') node[right] {$s_1'$};
\draw (s2') node[left] {$s_2'$};
\draw (r) node[right] {$r$};
\draw (s') node[above] {$s'$};
\draw (r') node[right] {$r'$};

\draw (v1) node[left] {$v_1$};
\draw (v2) node[right] {$v_2$};
\draw (v3) node[below] {$v_3$};
\node at (-1,.7) {$A_{s's_2'}$};
\node at (2.8,3.2) {$A_{rs}$};

\draw (v1') node[left] {$v_1'$};
\draw (v2') node[right] {$v_2'$};
\draw (w') node[above] {$w'$};
\draw (v3') node[right] {$v_3'$};
\node at (.3,2.1) {$P_1$};
\node at (3.3,1.2) {$P_2$};
\node at (1.3,1.2) {$P_3$};
\node at (2.8,1.8) {$Q$};

\end{tikzpicture}
\end{center}
\caption{Also for the proof of (8).} 
\label{fig:step9b}
\end{figure}

Now suppose that 
$s$ is a strict descendant of $s'$, and hence a descendant of one of $s_1',s_2'$ say $s_1'$. 
Let $P_1,P_2,P_3$ be three disjoint paths between the adhesions of
$s's_1'$ and $rs$, where $P_3$ joins $v_3', v_3$.  By exchanging $v_1,v_2$ if necessary we may assume that $P_1$ joins $w',v_1$ and $P_2$ joins $v_2',v_2$.
(See Figure \ref{fig:step9b}.) 
If $t$ exists, then the adhesion at $rt$ separates $A_{rt}$ from $B_{rt}$, and so, whether or not $t$ exists,
the interior of $Q$ is disjoint from $ A_{rs}, A_{s's_2'}$ and the interior of $P_3$.
By contracting all of $P_1,P_2,P_3, Q$, we convert the model of $\SP_d$
in $(G[A_{s's_2'}\setminus \{v_1'\}], v_3',w')$ and the model  of $\SP_d$
in $(G[A_{rs}\setminus \{v_1\}], v_2,v_3)$ to a model of $\SP_{d+1}$ in $(G[A_{r's'}],v_2',v_3')$, a contradiction. 
This proves (8).

\bigskip

For each $e=rs\in E(S)$, we define $J_e$ as follows.
If $s$ is a leaf of $S$, let $J_e$ be a triangle with vertex set the adhesion at $e$.
If $s$ is not a leaf of $S$, let $(v_1,v_2,v_3)$ be the standard numbering of the adhesion at $e$, and
let $J_e$ be the subgraph of $H(S)$ defined as follows. Let $X$ be the set of vertices of $S^-$ that are descendants of $s$. Start with the union of the graphs $G[W_{s'}]$ over all $s'\in X$,
and for every edge
$f$ of $T$ different from $e$ with exactly one end in $X$, add edges (if they are not already present) to make the adhesion at $f$ a clique. It follows that $J_e$ is a subgraph
of $H(S)$.

Let $e=rs\in E(S)$, and let $(v_1,v_2,v_3)$ be the standard numbering of the adhesion at $rs$. We say $e$ {\em behaves}
if: 
\begin{itemize}
\item there is a path $P_e$ of $J_e$ between $v_1,v_2$, and a tree $M_e$ of $J_e$ containing $v_3$, such that $V(P_e),V(M_e)$ are disjoint and have union $V(J_e)$;
\item every edge of $J_e$ either belongs to one of $P_e, M_e$, or has one end in $V(P_e)$ and the other in $V(M_e)$; and
\item there is a planar drawing of $J_e$ with $v_1,v_2,v_3$ all incident with the infinite region.
\end{itemize} 
We claim:
\\
\\
(9) {\em Every edge $rs$ of $S$ behaves.}
\\
\\
We prove that $e\in E(S)$ behaves, by induction on $|V(J_e)|$. Let $e = rs$, and let $W_s=\{v_1,v_2,v_3,w\}$, where $(v_1,v_2,v_3)$ is the standard numbering of the adhesion 
at $e$. If $s$ is a leaf of $S$, then from the definition of $J_e$, $w$ is not in $V(J_e)$, and $v_1,v_2,v_3$ are pairwise adjacent in $J_e$, and so $e$ behaves (taking $P_e$ to be
the path of length one between $v_1,v_2$, and $M_e$ the one-vertex tree with vertex $v_3$). 

Next, suppose that $s$ has exactly 
one child $s_1$ in $S$, and the adhesion of $ss_1$ is $W_s\setminus \{v_1\}$. Let $f = ss_1$. Then $v_1,w$ are adjacent in $H(S)$ and hence in $J_e$ since $G$ is three-connected. 
From the inductive hypothesis, $ss_1$ behaves,
so $P_{f}$ and $M_{f}$ exist, satisfying the bullets in the definition of ``behaves''. The ends of $P_f$ are $w,v_2$; let $P_e$ be obtained by adding $v_1$ and the edge $v_1w$ 
to $P_f$. Also, $v_3\in V(M_f)$; let $M_e=M_f$. The only edges of $J_e$ that are not edges of $J_f$ are $v_1w$ and possibly $v_1v_3$ or $v_3w$ (not $v_1v_2$ or $v_2w$, by (8)); and so $e$
behaves. 

Next, suppose that  $s$ has exactly 
one child $s_3$ in $S$, and the adhesion of $ss_3$ is $W_s\setminus \{v_3\}$. Let $f = ss_3$. Then $v_3,w$ are adjacent in $H(S)$ and hence in $J_e$ since $G$ is three-connected. 
From the inductive hypothesis, $ss_3$ behaves,
so $P_{f}$ and $M_{f}$ exist. The ends of $P_f$ are $v_1,v_2$; let $P_e=P_f$.
Also, $w\in V(M_f)$; let $M_e$ be obtained from $M_f$ by adding $v_3$ and the edge $v_3w$. The only edges of $J_e$ that are not edges of $J_f$ are between $v_3$ and 
$\{v_1,v_2,w\}$, and possibly $v_1w,v_2w$, so $e$
behaves.

Finally, suppose that $s$ has two children $s_1,s_2$, and the adhesions at $ss_i$ are $W_s\setminus \{v_i\}$ for $i =1,2$. Let $f_i = ss_i$ for $i = 1,2$.
From the inductive hypothesis,
$f_1$ and $f_2$ both behave, and so $P_{f_1},M_{f_1}, P_{f_2}, M_{f_2}$ exist and satisfy the bullets above. Since $P_{f_1}$ has ends $v_2,w$ and $P_{f_2}$
has ends $v_1,w$, their union
$P_{f_1}\cup P_{f_2}$ is a path between $v_1,v_2$; let $P_e$ be this path. Let $M_e = M_{f_1}\cup M_{f_2}$.  
Every edge of $J_e$ is an edge of one of $J_{f_1}$ or of $J_{f_2}$ except possibly $v_1v_3$ or $v_2v_3$  (because $v_1v_2, v_1w,v_2w$ are not new edges by (8)); and so $e$ behaves. This proves (9).

\bigskip

From (9), we deduce that $H(S)$ is near-outerplanar. This proves \ref{tw3full}.~\bbox
\section{The stopsign}
Let $H$ be a graph obtained from a cycle of length at least three and a path by identifying one end of the path with a vertex of the cycle. We call such a
graph a {\em stopsign}. Every stopsign is a subgraph of a fan, so we can use \ref{fanthm2} to give information about the graphs
that do not contain a given stopsign.
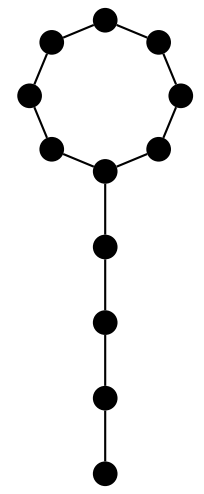
\begin{figure}[h!]
\begin{center}
\begin{tikzpicture}[scale = 1]
\tikzstyle{every node}=[inner sep=1.5pt, fill=black,circle,draw]
\def\t{45}
\def\r{.5}
\node (1) at (\r,0) {};
\node (2) at ({\r*cos(\t)}, {\r*sin(\t)}) {};
\node (3) at ({\r*cos(2*\t)}, {\r*sin(2*\t)}) {};
\node (4) at ({\r*cos(3*\t)}, {\r*sin(3*\t)}) {};
\node (5) at ({\r*cos(4*\t)}, {\r*sin(4*\t)}) {};
\node (6) at ({\r*cos(5*\t)}, {\r*sin(5*\t)}) {};
\node (7) at ({\r*cos(6*\t)}, {\r*sin(6*\t)}) {};
\node (8) at ({\r*cos(7*\t)}, {\r*sin(7*\t)}) {};

\node (9) at (0, {-\r-.5}) {};
\node (10) at (0, {-\r-1}) {};
\node (11) at (0, {-\r-1.5}) {}; 
\node (12) at (0, {-\r-2}) {};

\foreach \from/\to in {1/2,2/3,3/4,4/5,5/6,6/7,7/8,1/8,7/9,9/10,10/11,11/12}
\draw [-] (\from) -- (\to);

\end{tikzpicture}
\end{center}
\caption{A stopsign.} \label{fig:stopsign}
\end{figure}
We claim:
\begin{thm}\label{stopsign}
For every stopsign $H$, 
there exist $\ell, m$ such that every simple connected graph with no $H$ minor admits a rooted tree-decomposition
$\mac T=(T,(W_t:t\in V(T)))$ such that either:
\begin{itemize}
\item $\mac T$ has height at most $m$, and
one torso is a cycle, and all other torsos have at most $\ell+1$ vertices; or
\item $\mac T$ has width at most $\ell$, and every path of $T$ of length $m$ contains an edge with adhesion of size 1.
\end{itemize}
\end{thm}
\Proof 
Let $F$ be a fan containing
$H$ as a minor. Choose $\ell_1, m_1$ such that setting $\ell=\ell_1$ and $m=m_1$ satisfies \ref{fanthm2}. Choose $\ell_2,m_2$
such that every graph with no path of length $|V(F)|^2$ admits a tree-decomposition of height at most $m_2$ and width at most $\ell_2$.
Let $\ell =\max (\ell_1,\ell_2, |V(F)|)+1$ and $m = 2\max (m_1,m_2)+2$. We claim that $\ell, m$ satisfy \ref{stopsign}.

To see this, let $G$ be a simple connected graph with no $H$ minor.
Suppose first that $G$ has no cycle of length at least $|V(F)|$. Then each maximal two-connected component has no path of length $|V(F)|^2$ (by a theorem of Dirac~\cite{dirac}), 
and so admits a tree-decomposition of height at most $m_2$ and width at most $\ell_2$. By combining these tree-decompositions,
we deduce that $G$  admits a tree-decomposition of width at most $\ell_2$ such that every path of its tree of length $2m_2+1$ contains an edge with adhesion of size 1, and so the theorem holds.

Thus we may assume that some maximal two-connected component $G'$ has a cycle of length at least $|V(F)|$. For each $v\in V(G')$,
let $B_v$ be the component of $G\setminus \{v\}$ that contains $G'\setminus \{v\}$, and let $A_v=G\setminus V(B_v)$ (thus, $v\in V(A_v)$).
Hence the subgraphs $A_v\;(v\in V(G'))$ are pairwise vertex-disjoint, and their union, together with $G'$, equals $G$. Since $G$
does not contain $H$, it follows that for each $v\in V(G')$, there is no path in $A_v$ of length $|V(H)|$ with one end $v$,
and hence there is no path in $A_v$ of length $2|V(F)|$. Consequently, $A_v$ admits a tree-decomposition of height at most $m_2$ and width at most $\ell_2$.

Now, since $G'$ is two-connected and does not contain $F$ as a minor, \ref{fanthm2} implies that $G'$ admits a tree-decomposition $\mac T=(T,(W_t:t\in V(T)))$ of height at most $m_1$ such that for every $t\in V(T)$, either $|W_t|\le \ell_1+1$ or the torso of $\mac T$
at $t$ is a cycle.
If each $|W_t|\le \ell$, then by combining this tree-decomposition with that of $A_v$ for each $v\in V(G')$, we deduce that
$G$ admits a tree-decomposition of width at most $\ell$ such that every path of length more than $2\max (m_1,m_2)$
in $T$ contains an edge with adhesion of size one, and so the theorem holds. Thus, we assume that at least one torso of $\mac T$
is a cycle of  length at least $\ell\ge |V(F)|$. Since $G'$ is two-connected, there exists $X=\{x_1\LL x_n\}\subseteq V(G')$
with $n\ge \ell$, such that (reading subscripts modulo $n$):
\begin{itemize}
\item for each component $D$ of $G'\setminus X$, there exists $i\in \{1\LL n\}$ such that
$x_i, x_{i+1}$ have neighbours in $V(D)$ and no other members of $X$ have neighbours in $V(D)$; and
\item for $1\le i\le n$, either $x_ix_{i+1}$ is an edge of $G'$, or there is a component $D$ of $G'\setminus X$ such that $x_i, x_{i+1}$ have neighbours in $V(D)$.
\end{itemize}
Consequently, the same is true for $G$, except that there may be components $D$ of $G\setminus X$ such that only one vertex of $X$ has a neighbour in $V(D)$. Hence, for $1\le i\le n$ there is a subgraph $C_i$ of $G$ with the following properties,
again, reading subscripts modulo $n$ (see Figure \ref{fig:bangle}):
\begin{itemize}
\item $x_i, x_{i+1}\in V(C_i)$ and there is a path in $C_i$ between $x_i, x_{i+1}$;
\item for $1\le i<j\le n$, if $j\ne i+1$ and $(i,j)\ne (1,n)$ then $C_i\cap C_j$ is null;
\item $C_i\cap C_{i+1}=\{x_{i+1}\}$ for $1\le i\le n$; and
\item $C_1\cupcup C_n=G$.
\end{itemize}
\begin{figure}[h!]
\begin{center}
\begin{tikzpicture}[scale = 1.3]
\tikzstyle{every node}=[inner sep=1.5pt, fill=black,circle,draw]
\def\h{.8}
\def\t{30}
\node (1) at (1,0) {};
\node (2) at ({cos(\t)}, {sin(\t)}) {};
\node (3) at ({cos(2*\t)}, {sin(2*\t)}) {};
\node (4) at ({cos(3*\t)}, {sin(3*\t)}) {};
\node (5) at ({cos(4*\t)}, {sin(4*\t)}) {};
\node (6) at ({cos(5*\t)}, {sin(5*\t)}) {};
\node (7) at ({cos(6*\t)}, {sin(6*\t)}) {};
\node (8) at ({cos(7*\t)}, {sin(7*\t)}) {};
\node (9) at ({cos(8*\t)}, {sin(8*\t)}) {};

\draw[fill=yellow] (1) to [out= 10, in=\t-10, looseness = 5] (2);
\draw[fill=yellow] (2) to [out= \t+10, in=2*\t-10, looseness = 5] (3);
\draw[fill=yellow] (3) to [out= 2*\t+10, in=3*\t-10, looseness = 5] (4);
\draw[fill=yellow] (4) to [out= 3*\t + 10, in=4*\t-10, looseness = 5] (5);
\draw[fill=yellow] (5) to [out= 4*\t+10, in=5*\t-10, looseness = 5] (6);
\draw[fill=yellow] (6) to [out= 5*\t+10, in=6*\t-10, looseness = 5] (7);
\draw[fill=yellow] (7) to [out= 6*\t+10, in=7*\t-10, looseness = 5] (8);
\draw[fill=yellow] (8) to [out= 7*\t+10, in=8*\t-10, looseness = 5] (9);

\draw[fill=white] (1) to [out= 70, in=\t-70, looseness = 1.5] (2);
\draw[fill=white] (2) to [out= \t+70, in=2*\t-70, looseness = 1.5](3);
\draw[fill=white] (3) to [out= 2*\t+70, in=3*\t-70, looseness = 1.5](4);
\draw[fill=white] (4) to [out= 3*\t+70, in=4*\t-70, looseness = 1.5](5);
\draw[fill=white] (5) to [out= 4*\t+70, in=5*\t-70, looseness = 1.5](6);
\draw[fill=white] (6) to [out= 5*\t+70, in=6*\t-70, looseness = 1.5](7);
\draw[fill=white] (7) to [out= 6*\t+70, in=7*\t-70, looseness = 1.5](8);
\draw[fill=white] (8) to [out= 7*\t+70, in=8*\t-70, looseness = 1.5](9);

\def\r{1.2}
\draw[domain=245:355,smooth,variable=\x,thick, dotted] plot ({\r*cos(\x)},{\r*sin(\x)});

\tikzstyle{every node}=[]
\node at (.3,1.2) {$C_i$};
\draw (3) node[below] {$x_i$};
\draw (4) node[below] {$x_{i+1}$};
\end{tikzpicture}
\end{center}
\caption{For \ref{stopsign}.} \label{fig:bangle}
\end{figure}
Choose $n$ maximum such that there is such a set $X=\{x_1\LL x_n\}$ and $C_1\LL C_n$. It follows that for $1\le i\le n$, either $x_i,x_{i+1}$ are adjacent, or there are two paths of $C_i$ between $x_i, x_{i+1}$ with disjoint interiors. We claim that for $1\le i\le n$,
$C_i$ has no path of length $|V(F)|^2$. Suppose that $i=1$ say, and there is a path $P$ in $C_1$ of length $|V(F)|^2$.
There is a path in $C_1$ between $x_1,x_2$, say $Q$. There is no path $R$ of $C_1$ with $|V(R)|\ge |V(F)|$ with one end in $Q$ and no
other vertex in $Q$, since $G$ does not contain $H$ as a minor. (Let us call this fact the ``forbidden path condition''.)
If $P\cap Q$ is null, then the union of $P$ and a path of $C_1$ between $P,Q$
gives a path that violates the forbidden path condition, a contradiction.
Thus $P\cap Q$ is non-null. If any component of $P\setminus V(Q)$ has at least $|V(F)|-1$ vertices,
then again the forbidden path condition is violated.
So there is no such component. Let $|V(P\cap Q)|=k$; then $P\setminus V(Q)$ has at most $k+1$ components, each with at most $|V(F)|-2$ vertices, and
so 
$$|V(F)|^2 \le |V(P)|\le k+(k+1)(|V(F)|-2)< (k+1)(|V(F)|-1).$$
It follows that $k>|V(F)|$, and 
in particular $|V(Q)|> |V(F)|$.
Hence every path of $C_1$ between $x_1,x_2$ has at least $|V(F)|+1$ vertices. So $x_1,x_2$ are nonadjacent, and there are two paths of $C_1$
between $x_1,x_2$ with disjoint interiors, both of length at least  $|V(F)|$, and again the forbidden path condition is violated.
This proves that $C_i$ has no path of length $|V(F)|^2$, for $1\le i\le n$. Hence $C_i$ admits a tree-decomposition
of height at most $m_2$ and width at most $\ell_2$, for $1\le i\le n$; and so it admits a tree-decomposition
of height at most $m_2$ and width at most $\ell_2+1$, such that some bag contains both $x_i, x_{i+1}$. Consequently, $G$ admits a  rooted tree-decomposition
of height at most $m_2+1$ such that the torso at the root is a cycle and all others have at most $\ell_2+2$ vertices, and the theorem holds.
This proves \ref{stopsign}.~\bbox
\section{Excluding caterpillars}

A {\em caterpillar} is a tree such that some path contains a neighbour of every vertex not in the path. 
We know a rough structure theorem for
excluding a tree, namely \ref{GM1}, but what if the tree we are excluding is a caterpillar? Is there a rough structure theorem for that? Just saying that path-width is
bounded is not enough: it excludes no caterpillars, because all caterpillars have path-width at most one. But there is a rough structure theorem that we prove
in this section. We first need some definitions. 

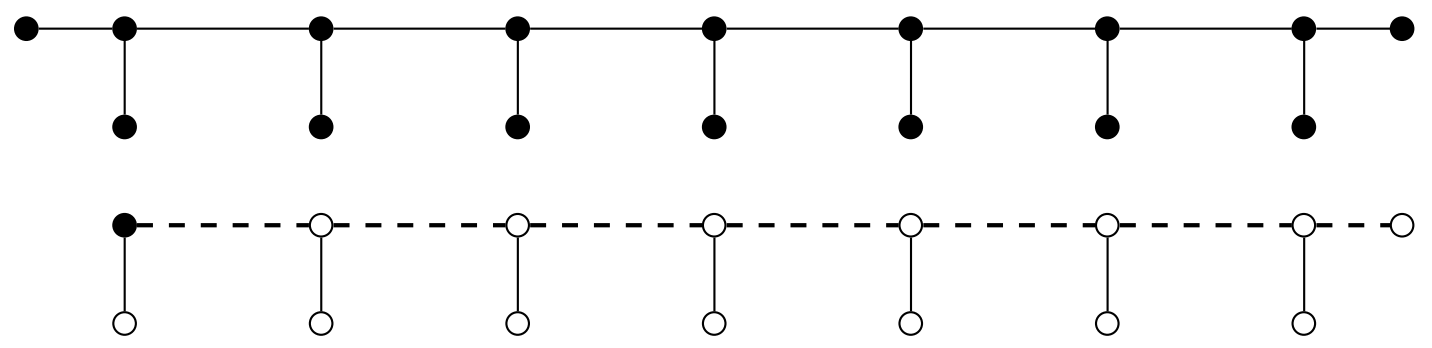
\begin{figure}[h!]
\begin{center}
\begin{tikzpicture}[scale = 1.3]
\tikzstyle{every node}=[inner sep=1.5pt, fill=black,circle,draw]
\node (1) at (1.5,1) {};
\node (2) at (2,1) {};
\node (3) at (3,1) {};
\node (4) at (4,1) {};
\node (5) at (5,1) {};
\node (6) at (6,1) {};
\node (7) at (7,1) {};
\node (8) at (8,1) {};
\node (9) at (8.5,1) {};
\node (2') at (2,.5) {};
\node (3') at (3,.5) {};
\node (4') at (4,.5) {};
\node (5') at (5,.5) {};
\node (6') at (6,.5) {};
\node (7') at (7,.5) {};
\node (8') at (8,.5) {};

\foreach \from/\to in {1/2,2/3,3/4,4/5,5/6,6/7,7/8,8/9,2/2',3/3',4/4',5/5',6/6',7/7',8/8'}
\draw [-] (\from) -- (\to);

\node (2) at (2,0) {};
\tikzstyle{every node}=[inner sep=1.5pt, fill=white,circle,draw]
\node (3) at (3,0) {};
\node (4) at (4,0) {};
\node (5) at (5,0) {};
\node (6) at (6,0) {};
\node (7) at (7,0) {};
\node (8) at (8,0) {};
\node (9) at (8.5,0) {};
\node (2') at (2,-.5) {};
\node (3') at (3,-.5) {};
\node (4') at (4,-.5) {};
\node (5') at (5,-.5) {};
\node (6') at (6,-.5) {};
\node (7') at (7,-.5) {};
\node (8') at (8,-.5) {};

\foreach \from/\to in {2/2',3/3',4/4',5/5',6/6',7/7',8/8'}
\draw [-] (\from) -- (\to);

\foreach \from/\to in {2/3,3/4,4/5,5/6,6/7,7/8,8/9}
\draw [thick, dashed] (\from) -- (\to);

\end{tikzpicture}
\end{center}
\caption{The caterpillar $L_7$, and a one-rooted graph of type $7$.} \label{fig:caterpillar}
\end{figure}

For $d\ge 1$, let $L_d$ be the caterpillar formed by taking a path with $d+2$ vertices and adding $d$ new vertices, one adjacent to each internal vertex of the path.
(See Figure \ref{fig:caterpillar}.)
Every caterpillar is a minor of $L_d$ for some $d$, so for our purposes it suffices to exclude $L_d$. 

For $d\ge 1$, we say a one-rooted graph $(H,v)$ has {\em type $d$} if $H$ is a tree and can be obtained from a path with ends $v,v'$ by adding $d$ new vertices, all  
adjacent to distinct vertices of the path, and not adjacent to $v'$.

A {\em tree-partition} of a graph $G$ is a pair $(T,(W_t:t\in V(T)))$, such that 
\begin{itemize}
\item $T$ is a rooted tree;
\item the sets $W_t\;(t\in V(T))$ are pairwise disjoint subsets of $V(G)$ with union $V(G)$;
\item for each edge $e=uv$ of $G$, let $s,t\in V(T)$ with $u\in W_s$ and $v\in W_t$; then one of $s,t$ is an ancestor of the other.
\end{itemize}
Its {\em height} is the height of $T$.

We claim:
\begin{thm}\label{catthm}
For every $d\ge 1$, there exists $\ell_d$ such that every simple connected graph with no $L_d$ minor admits a tree-partition $(T, (W_t:t\in V(T)))$ of 
height at most $d-1$, 
such that for each $t\in V(T)$, if the distance in $T$ between $t$ and the root is less than $d-1$ then $|W_t|\le \ell_d$, and if this distance is $d-1$ then $G[W_t]$
is a path or cycle.
\end{thm}
\Proof Choose $k$ such that every graph with no $L_d$ minor has tree-width less than $k$. 
We will prove, by induction on $h$, that for all $h$ with $1\le h\le d$, 
every simple connected graph with no $L_h$ minor admits a tree-partition $(T, (W_t:t\in V(T)))$ of
height at most $h-1$, 
such that for each $t\in V(T)$, if the distance in $T$ between $t$ and the root is less than $h-1$ then $|W_t|\le k$, and if this distance is $h-1$ then $G[W_t]$
is a path or cycle. Consequently, this will show that setting $\ell_d = k$ satisfies the theorem.

Every simple connected graph with no $L_1$ minor is a path or cycle, so the result holds when $h=1$, taking $|V(T)|=1$.
We assume that $h\ge 2$, and the claim holds for $h-1$.  Now $G$ contains $L_h$ as a minor if and only if it has a subgraph that is a subdivision of $L_h$,
since $L_h$ has maximum degree three; and the same for $L_{h-1}$. 
\\
\\
(1) {\em Either there are two vertex-disjoint subdivisions of $L_{h-1}$ in $G$,
or there exists $X\subseteq V(G)$ with $|X|\le k$ such that there is no subdivision of $L_{h-1}$ in $G\setminus X$.}
\\
\\
This is true by the 
Erd\H{o}s-P\'osa property of planar minors proved in~\cite{GM5}, but the proof is short so we give it anyway. Take a tree-decomposition of $G$ of width less than $k$.
For each subdivision $S$ of $L_{h-1}$ in $G$, let $T_S$ be the subtree induced on the set of $t\in V(T)$ with $W_t\cap V(S)\ne \emptyset$. If there exist 
subdivisions $S,S'$ of $L_{h-1}$ such that $T_S,T_{S'}$ are vertex-disjoint, then $S,S'$ are vertex-disjoint. If not, then by the Helly property of subtrees of a tree,
there exists $t\in V(T)$ such that $t\in V(T_S)$ for every subdivision $S$ of $L_{h-1}$ in $G$; and then there is no subdivision of $L_{h-1}$ in $G\setminus W_t$.
This proves (1). 

\bigskip

We also need:
\\
\\
(2) {\em Let $S$ be a subdivision of $L_{h-1}$; then for every $u\in V(S)$, the one-rooted graph $(S,u)$ has a one-rooted subgraph of type $\lceil h/2\rceil$.}
\\
\\
The proof is easy case-checking (using that $h\ge 2$) and we leave it to the reader.
\bigskip

Now for the inductive proof of the claim: let $G$ be 
a connected simple graph with no $L_h$ minor. Suppose that there are two vertex-disjoint subdivisions $S,S'$ of $L_{h-1}$ in $G$, and let $P$ be a minimal path
in $G$ between them, with ends $u\in V(S)$ and $u'\in V(S')$. By (2), $(S,u)$ has a one-rooted subgraph of type $\lceil h/2\rceil$, and the same for
$(S',u')$; and joining these by $P$ gives a subgraph that contains a subdivision of $L_h$, a contradiction. So there do not exist
two vertex-disjoint subdivisions of $L_{h-1}$ in $G$. By (1), 
there exists $X\subseteq V(G)$ with $|X|\le k$ such that there is no subdivision of $L_{h-1}$ in $G\setminus X$, and hence $G\setminus X$ has no $L_{h-1}$ minor.
From the inductive hypothesis,
each component $C$ of $G\setminus X$ admits a rooted tree-partition $(T_C,(W_t:t\in V(T_C)))$ of height at most $h-2$, such that for each $t\in V(T_C)$,
if the distance in $T_C$ between $t$ and the root is less than $h-2$ then $|W_t|\le k$, and if this distance is $h-2$ then $G[W_t]$
is a path or cycle.
Add a new vertex $t_0$ adjacent to the root
of $T_C$ for each $C$, making a tree $T$; and define $W_{t_0} = X$. Then $(T,(W_t:t\in V(T)))$ satisfies the theorem. This proves \ref{catthm}.~\bbox
\section{The converses}\label{sec:converses}

As promised earlier, in this section we sketch the proof of the ``rough converses'' to \ref{SPthm}, \ref{outerpthm}, \ref{innerpthm}, \ref{tw3full}, \ref{fanthm}, \ref{stopsign} and \ref{catthm}. 
Let $Q$ be a connected simple graph; we define graphs of {\em type $(Q,a,k)$} where $a\ge 1$ and $k\ge 0$ as follows. Graphs of type $(Q,a,0)$ are just copies of $Q$. For $k\ge 1$, a graph $G$ is of 
type $(Q,a,k)$ if it can be constructed from a copy $J$ of $Q$ as follows. 
For each vertex $v\in V(J)$, take $a+1$ pairwise vertex-disjoint graphs each of type $(Q,a,k-1)$, choose a vertex from each of them, and identify these $a+1$ selected 
vertices with $v$.

We need the following lemma:
\begin{thm}\label{AIlemma}
Let $Q,a,k$ be as above, and let $H$ be a graph of type $(Q,a,k)$. Let $G$ be a graph that admits a tree-decomposition of height at most $k$, with all adhesions of size at most $a$, and such that no torso contains $Q$ as a minor. Then $G$ does not contain $H$ as a minor.
\end{thm}
\Proof We proceed by induction on $k$. The claim for $k=0$ is clear, so we assume that $k\ge 1$ and the result holds for $k-1$. Let 
$G$ be a graph that admits a tree-decomposition $\mac T = (T,(W_t:t\in V(T)))$ of height at most $k$ and with all adhesions of size at most $a$,
such that no torso contains $Q$
as a minor. Suppose for a contradiction that there is a model $\eta$ of $H$ in $G$. 

Now $H$ is obtained from a copy $J$ of $Q$, by attaching $a+1$ graphs of type $(Q,a,k-1)$ to each vertex $v$ of $J$. Let $r_0$ be the root of $T$, and suppose that 
$V(\eta(v))\cap W_{r_0}\ne \emptyset$ for each $v\in V(J)$. Let $L$ be the torso of $\mac T$ at $r_0$, and let $\eta'(v) = L[V(\eta(v))\cap W_{r_0}]$
for each $v\in V(J)$. Hence for each $v\in V(J)$, the graph $\eta'(v)$ is a non-null connected subgraph of $L$.
Moreover, for each edge $uv$ of $J$, we claim that there is an edge of $L$ between $V(\eta'(u))$ and $V(\eta'(v))$. To show this, certainly there is an edge $xy$
of $G$ between  some $x\in V(\eta(u))$ and some $y\in V(\eta(v))$. If both $x,y\in W_{r_0}$ then the claim holds, so we assume that $x\notin W_{r_0}$; and hence $x$
belongs to $A_f\setminus W_{r_0}$ for some edge $f$ of $T$ incident with $r_0$. Consequently $x,y$ both belong to $V(a_f)$, and so  $V(\eta'(u))$ and $V(\eta'(v))$
both have a vertex in the adhesion of $f$. Since the adhesion of $f$ is a clique of $L$, this proves our claim that  there is an edge of $L$ between $V(\eta'(u))$ and $V(\eta'(v))$. Consequently, $L$ contains $Q$ as a minor, a contradiction. 

Thus, there exists $v\in V(J)$ such that $V(\eta(v))\cap W_{r_0}=\emptyset$, and therefore $V(\eta(v))\subseteq A_f)\setminus W_{r_0}$
 for some edge $f$ of $T$ incident with $r_0$. In $H$ there are $a+1$ subgraphs $J_1\LL J_{a+1}$,  each of type $(Q,a,k-1)$,  each containing $v$, and otherwise 
pairwise vertex-disjoint. For $1\le i\le a+1$, let $V_i=\bigcup_{z\in V(J_i)}V(\eta(z))$. Thus, $G[V_i]$ is connected for each $i$, and the sets
$V_1\LL V_{a+1}$ pairwise have intersection $V(\eta(v))$. 

Since the adhesion ($C$, say) at $f$ in $\mac T$ has size at most $a$, and is disjoint from $V(\eta(v))$, there exists $i\in \{1\LL a+1\}$ such that 
$V_i\cap C=\emptyset$; and since $G[V_i]$ is connected, it follows that $V_i\subseteq A_f\setminus W_{r_0}$. But $G[V_i]$ contains $J_i$ as a minor. Hence so does 
$G[A_f]$; but the latter admits a tree-decomposition of height at most $k-1$, with all adhesions of size at most $a$, and such that no torso contains $Q$ as a minor,
contradicting the inductive hypothesis. This proves \ref{AIlemma}.~\bbox

For \ref{SPthm}, we need to show:
\begin{thm}\label{SPconv}
For all $m\ge 0$ there is a series-parallel graph $H$, such that no simple two-connected graph with an $H$ minor admits a rooted tree-decomposition with adhesion two,
of height at most $m$, and such that each torso is either outerplanar or innerplanar.
\end{thm}
\Proof Let $Q$ be the graph $\SP_5$ (with roots ignored); so $Q$ is connected, and neither outerplanar nor innerplanar. Let $H$ be a graph of type $(Q,3,m)$; then $H$ is series-parallel. 
By \ref{AIlemma}, if $G$ admits a rooted tree-decomposition with adhesion two,
of height at most $m$, and such that each torso is either outerplanar or innerplanar, then $H$ is not a minor of $G$. This proves \ref{SPconv}.~\bbox

What about the other converses? For \ref{outerpthm}, the same method works, taking $Q$ to be a connected outerplanar graph that is not innerplanar, 
such as $\OP_7$, or the graph obtained from a cycle of length six by making its first, third and fifth vertices pairwise adjacent. (Then graphs of type $(Q,3,m)$ are 
outerplanar, and not minors of any graph $G$ that admits a decomposition as in \ref{outerpthm}, by \ref{AIlemma}.)
For \ref{innerpthm}, again the same method works, taking $Q$ to be a connected innerplanar graph that is neither outerplanar nor an apex-forest, such as $\IP_5$.
(Then graphs of type $(Q,3,m)$ are also innerplanar, and not a minor of any graph $G$ that admits a decomposition as in \ref{innerpthm}.) For \ref{tw3full}, take $Q$ 
to be a connected series-parallel graph that is not of any of the four types listed in \ref{tw3full}. For instance:
\begin{itemize}
\item  let $Q_1$ be a connected series-parallel graph that is not near-outerplanar, such as the graph obtained from $K_{2,3}$ by identifying each of its three vertices of degree two with a vertex of some new triangle;
\item let $Q_2$ be a connected series-parallel graph obtained from two disjoint copies of $\SP_7$ by adding an edge between them (so $Q_2\setminus \{v\}$ is neither 
outerplanar nor innerplanar, for each vertex $v$); and
\item let $Q_3$ be a tree with path-width at least four.
\end{itemize}
Let $Q$ be obtained from the disjoint union of $Q_1,Q_2,Q_3$ by choosing one vertex from each component and identifying them. 
(Then graphs of type $(Q,3,m)$ are also series-parallel, and not a minor of any graph $G$ that admits a decomposition as in \ref{tw3full}.)

For a converse to \ref{fanthm}, it suffices to show:
\begin{thm}\label{fanthmconv}
For all $\ell, m\ge 1$ there is a fan $H$ such that no simple two-connected graph $G$ with an $H$ minor admits
a tree-decomposition $(T,(W_t:t\in V(T)))$ of height at most $m$ such that for each $t\in V(T)$, either $|W_t|\le \ell$ or the torso at $t$ is a cycle.
\end{thm}
\noindent{\bf Proof sketch.\ \ }
To see this, we may assume that $\ell\ge 3$ by increasing $\ell$ if necessary. Choose $n$ such that no
graph with a path of length $n$ admits a tree-decomposition of width at most $\ell$ and height at most $m$. Let $H$ be a fan of
length $n$. 
It suffices to show that $H$ does not admit a decomposition as in the theorem.
Suppose that it does,  that is, 
it admits a tree-decomposition $\mac T=(T,(W_t:t\in V(T)))$ of height at most $m$ such that for each $t\in V(T)$, either $|W_t|\le \ell$ or the torso at $t$ is a cycle.
Since $H$ is a chordal graph, any torso that is a cycle has length three, and hence
has three vertices, and hence at most $\ell$ vertices. Thus $\mac T$ has width at most $\ell$, and so $H$ has no path of length $n$, a contradiction. This proves \ref{fanthmconv}.~\bbox

For \ref{stopsign}, we show:
\begin{thm}\label{stopsignconv}
For all $\ell, m\ge 1$ there is a stopsign $H$ such that no simple connected graph $G$ with an $H$ minor admits
a rooted tree-decomposition $\mac T=(T,(W_t:t\in V(T)))$ such that either
\begin{itemize}
\item $\mac T$ has height at most $m$, and
one torso is a cycle, and all other torsos have at most $\ell+1$ vertices; or
\item $\mac T$ has width at most $\ell$, and every path of $T$ of length $m$ contains an edge with adhesion of size 1.
\end{itemize}
\end{thm}
\noindent{\bf Proof sketch.\ \ }
Choose $t$ such that no graph with a path of length $t$
admits a tree-decomposition of width at most $\ell$ and height at most $m$. Let  $H$ be a stopsign made
by identifying one end of a path of length $t+1$ with a vertex of a cycle of length $t+1$. 
As usual, it suffices to show that $H$ does not admit a tree-decomposition as in the theorem. Suppose 
$H$ admits a rooted tree-decomposition $\mac T$ such that either:
\begin{itemize}
\item $\mac T$ has height at most $m$, and
one torso is a cycle, and all other torsos have size at most $\ell+1$; or
\item $\mac T$ has width at most $\ell$, and every path of $T$ of length $m$ contains an edge with adhesion of size 1.
\end{itemize}
In the second case, all cycles have bounded length, so if $t$ is sufficiently large then $H$ admits no such decomposition. In
the first case, $H$ admits a decomposition as shown in Figure \ref{fig:bangle}: there are subgraphs $C_1\LL C_n$ for some $n\ge 3$,
arranged in a circular order,
each sharing one vertex with the next and otherwise disjoint, with union $H$, and every path in each $C_i$ has length less than $t$.
In particular, no $C_i$ has a cycle of length at least $t+1$, so the (unique) cycle of $H$ contains all the
vertices that belong to more than one of the $C_i$'s. But there is a path of $H$ of length $t$ disjoint from the cycle of $H$, so it is a subgraph of some $C_i$, which is impossible. This proves \ref{stopsignconv}. ~\bbox

For a converse to \ref{catthm}, we need:
\begin{thm}\label{catthmconv}
For every $\ell,m\ge 0$, let $d(\ell,m) = 3(2(\ell+1))^m$; then no simple connected graph with an $L_{d(\ell,m)}$ minor admits a tree-partition $(T, (W_t:t\in V(T)))$ of
height at most $m$,
such that for each $t\in V(T)$, if the distance in $T$ between $t$ and the root is less than $m$ then $|W_t|\le \ell$, and if this distance is $m$ then $G[W_t]$
is a path or cycle.
\end{thm}
\noindent{\bf Proof sketch.\ \ }
It suffices to show that $L_{d(\ell,m)}$ does not admit such a tree-partition. This is true if $m=0$, and if $m\ge 1$ it follows immediately from an inductive hypothesis
on $m$, by deleting the bag of the root of $T$, using that if $X$ is a set of vertices of size at most $\ell$, then $L_{d(\ell,m)}\setminus X$
contains a copy of $L_{d(\ell,m-1)}$. This proves \ref{catthmconv}.~\bbox
\section*{AI use}
There was no use of AI in carrying out the research described in this paper. After writing it, we used AI for proof-reading, and AI suggested an improvement in a proof
in the final section that was so nice we could not resist including it. 


\begin{thebibliography}{99}

\bibitem{chuzhoy} J. Chuzhoy and Z. Tan, ``Towards tight(er) bounds for the excluded grid theorem,''
{\em J. Combinatorial Theory, Ser. B}, {\bf 146} (2021), 219--265.

\bibitem{diestel2} R. Diestel,  ``Graph minors I: a short proof of the path-width theorem'', {\em Combinatorics,
Probability and Computing}, {\bf 4} (1995), 27--30.


\bibitem{diestel} R. Diestel, T. Jensen, K. Gorbunov, and C. Thomassen, ``Highly connected sets and the excluded grid
theorem'', {\em J. Combinatorial Theory, Ser. B}, {\bf 75} (1999), 61--73.


\bibitem{wheel} G. Ding, S. Dziobiak, and H. Wu, ``Large $W_k$- or $K_{3,t}$-minors in 3-connected graphs'', {\em J. Graph Theory}, {\bf 82} (2016), 207--217.

\bibitem{dirac} G. Dirac, ``Some theorems on abstract graphs'', {\em Proc. London Math. Soc.} {\bf (3)2} (1952), 69--81.


\bibitem{erdos} P. Erd\H{o}s and L. P\'osa,  ``On independent circuits contained in a graph'', {\em Canadian Journal of Mathematics} {\bf 17} (1965), 347--352. 

\bibitem{haj} M. Hajiaghayi, ``A polynomial-time sixth-power excluded grid theorem'', {\tt arXiv:2611.16755}.



\bibitem{GM1} N. Robertson and P. Seymour,  ``Graph minors. I. Excluding a forest'',
{\em J. Combinatorial Theory, Ser. B,} {\bf 35} (1983), 39--61.



\bibitem{GM5} N. Robertson and P. Seymour,  ``Graph minors. V. Excluding a planar graph'',
{\em J. Combinatorial Theory, Ser. B,} {\bf 41} (1986), 92--114.


\bibitem{GM16} N. Robertson and P. Seymour, ``Graph minors. XVI. Excluding a non-planar graph'',
{\em J. Combinatorial Theory, Ser. B}, {\bf 89} (2003), 43--76.



\bibitem{sengupta} R. Sengupta, {\em Hedgewidth In Series-Parallel Graphs}, Senior Thesis, Princeton, 2012.



\bibitem{whitney} H. Whitney, ``2-isomorphic graphs'', {\em American J. Math.}, {\bf 55} (1933), 245--254.

\end{thebibliography}
\end{document}